\documentclass[final,onefignum,onetabnum]{siuro250211}

\usepackage{amsfonts}
\usepackage{graphicx}
\usepackage{epstopdf}
\usepackage{algorithmic}
\usepackage{mathtools}
\usepackage{amsfonts, amssymb}
\usepackage{float}
\usepackage{placeins}
\ifpdf
  \DeclareGraphicsExtensions{.eps,.pdf,.png,.jpg}
\else
  \DeclareGraphicsExtensions{.eps}
\fi

\usepackage{bm}
\newcommand{\vect}[1]{\ensuremath{\bm{#1}}}

\usepackage{enumitem}
\setlist[enumerate]{leftmargin=.5in}
\setlist[itemize]{leftmargin=.5in}

\newsiamremark{remark}{Remark}
\newsiamremark{hypothesis}{Hypothesis}
\crefname{hypothesis}{Hypothesis}{Hypotheses}
\newsiamthm{claim}{Claim}
\newsiamremark{fact}{Fact}
\crefname{fact}{Fact}{Facts}

\headers{Low-Rank Acceleration of the Operator Fourier Transform}{J. Kelley}

\title{Low-Rank Acceleration of the Operator Fourier Transform}

\author{Jack Kelley\thanks{Computational Modeling and Data Analytics Program, Virginia Tech
  (\email{jackkelley@vt.edu}).}}

\dedication{\small\textit{Project advisor: Daniel Appel\"o\thanks{Department of Mathematics, Virginia Tech (\email{appelo@vt.edu}).}}}

\usepackage{amsopn}

\ifpdf
\hypersetup{
  pdftitle={Low-Rank Acceleration of the Operator Fourier Transform},
  pdfauthor={J. Kelley}
}
\fi

\begin{document}

\maketitle

\begin{abstract}
We develop a numerical algorithm for the efficient solution or approximation of solutions to the Helmholtz equation on a structured grid in two dimensions. We make use of the Operator Fourier Transform (OFT) and a low-rank cross approximation scheme (Cross-DEIM) to decompose the problem into an integral over a pseudo-time of solutions to the Schr\"odinger equation. The OFT is a framework for solving operator equations like fractional Laplacian equations or the Helmholtz equation, when the latter is written as a product of two paraxial operators. The main computational cost in
the OFT is the solution to the Schr\"{o}dinger equation, especially when the dimension or mesh resolution is high. In this work, we alleviate this cost by utilizing a low-rank method. Such methods aim to beat the curse of dimensionality when low-rank structures are present in the solution. For a problem with a constant coefficient we show that the resulting solver converges at second order in the pseudo-time step and reaches the accuracy of the direct solution of the discretized linear system, which is limited by the second-order finite-difference discretization in space, while keeping the solution in low-rank form. The ranks of the accumulated solution change little as the grid is refined, which suggests a cost per time step that grows nearly linearly with the number of grid points in each direction.
\end{abstract}

\section{Introduction}
The study of waves arises in many science and engineering applications. Waves are able to carry information over a distance, enabling both natural (through animal sonar and speech) and technological communication \cite{LR-WaveHoltz}. Waves also underpin various imaging and probing methods, either through direct reception (eyes, cameras) or by solving an inverse problem as in, for example, computed tomography \cite{WaveHoltz1}. It is therefore necessary to develop fast and robust methods for computing the solution to problems involving waves. In this paper we are interested in solving the Helmholtz equation, which is the time-harmonic form of the wave equation.

We present an adaptation of the Operator Fourier Transform (OFT) \cite{Cubillos2024OFT} that makes use of the Cross-DEIM algorithm proposed in \cite{appelo2025lraa} for the low-rank solution of the Helmholtz equation. The specific formulation of the Helmholtz equation that we will use is
\begin{equation}
    \label{eq:helmholtz}
    \begin{aligned}
        m\left(\vect{x}\right)v\left(\vect{x}\right) + \dfrac{1}{\kappa^2}\Delta v\left(\vect{x}\right) &= g\left(\vect{x}\right), && \vect{x} \in \Omega, \\[4pt]
        v\left(\vect{x}\right) &= 0, && \vect{x} \in \partial\Omega,
    \end{aligned}
\end{equation}
where $\kappa$ is the free-space wavenumber, $m$ is the square of the index of refraction, and $\vect{x} \in \Omega \subset \mathbb{R}^d$ with $d$ the dimension. The Laplacian of $v\left(\vect{x}\right)$ is denoted by $\Delta v\left(\vect{x}\right)$. Classical iterative methods are known to converge slowly, or not at all, for the Helmholtz equation at high wavenumbers, since the discretized operator is indefinite \cite{ernst2012difficult}. Considerable progress has been made with specialized preconditioners, such as sweeping preconditioners \cite{engquist2011sweeping}, the method of polarized traces \cite{zepeda2016polarized}, and the related domain decomposition methods reviewed in \cite{gander2019class}. The OFT offers a different approach, as it leads to the implementation of a direct solver where the Helmholtz operator is inverted via the solution of two time-dependent linear Schr\"odinger equations. These solutions are then summed within a Fourier integral over a time-like variable $\tau$.

The OFT framework allows for any method to be used for the evaluation of the Schr\"odinger equation. We choose to use low-rank methods to compute an approximation to the solution of the Schr\"odinger equation. As described in \cite{lubich2007dynamiclowrank}, the Schr\"odinger equation has a long history of being evaluated on low-rank manifolds, stretching back to Dirac in the 1930s \cite{dirac1930exchange}. It is demonstrated in Section II of \cite{lubich2008quantumbook} that for the Schr\"odinger equation, when a large system is approximated by a low-rank representation, the error is bounded by the best-approximation error of the Schr\"odinger solution on to the low-rank manifold.

The goal of low-rank methods is to exploit a low-rank structure present in the solution. On a grid in two dimensions the discrete solution can be represented by the matrix $X \in \mathbb{C}^{m\times n}$ where the entry $X_{ij}$ approximates the solution at the grid point $\left(x_i, y_j\right)$. We write $X = V + iW$ with $V, W \in \mathbb{R}^{m \times n}$, as detailed in \cref{sec:low-rank}, and approximate both $V$ and $W$ as a rank $r$ matrix, where $r \ll \min\left(m, n\right)$. This allows us to only store and evaluate part of the solution, reducing both the memory footprint and the number of operations required for a solution. Low-rank acceleration is a natural extension of the OFT, as the primary computational cost of the evaluation of the OFT integral is the solution to the Schr\"odinger equation. Using the Cross-DEIM algorithm, the cost of one time step of our method on an $n \times n$ grid is $O\left(n r \log n + n r^2\right)$, compared with $O\left(n^2 \log n\right)$ when full matrices are used; see \cref{sec:lr_oft}.

To advance $V$ and $W$ in time we introduce a new scheme making use of this real-imaginary splitting. We employ an exponential time-stepping scheme developed in \cite{Cubillos2024OFT} in which the initial and final step sizes are prescribed and the intermediate nodes grow geometrically; this is well-suited to the OFT integrand, whose decay in $\tau$ allows larger steps at later times without significant loss of accuracy. At each step we apply a midpoint discretization and diagonalize the discrete Laplacian via a discrete sine transform, and in this basis, for a constant coefficient $m$, the coupled update for $\hat{V}$ and $\hat{W}$ decouples pointwise. The time step then reduces, at each grid point $(i, j)$, to multiplication by a $2\times 2$ rotation matrix with angle $\theta$ and prefactor $\phi$ determined by the coefficient $m$ and the eigenvalues of the discrete Laplacian. When the medium is non-absorbing, $\phi=1$ and the scheme is norm-preserving; otherwise $\phi < 1$ leads to dissipation.

The low-rank algorithm we choose to use is Cross-DEIM, as developed in \cite{appelo2025lraa}. We choose to use this algorithm as it does not require us to ever form the full matrices used in our Schr\"odinger solver, keeping all intermediate steps low-rank. The Cross-DEIM algorithm is an iterative greedy method that finds a low rank approximation in SVD form of a matrix $G$, without requiring access to the full matrix. Given a set of row indices $\mathcal{I}$ and column indices $\mathcal{J}$ of the matrix $G$, we approximate $G$ as
\begin{equation}
    G \approx G(:, \mathcal{J})\, M\, G(\mathcal{I},:), \qquad M = G(\mathcal{I}, \mathcal{J})^\dagger,
\end{equation}
where $G(:, \mathcal{J})$ denotes the submatrix of $G$ formed by all rows and the columns indexed by $\mathcal{J}$, and $G(\mathcal{I}, :)$ denotes the rows indexed by $\mathcal{I}$ and all columns. Unlike the related CUR decomposition \cite{mahoney2009cur}, this approximation does not require access to the full matrix $G$. The selection of $\mathcal{I}$ and $\mathcal{J}$ and the stable computation of the approximation are described in \cref{sec:cross-deim}. One advantage of using Cross-DEIM is that it is rank adaptive, i.e. it is possible to prescribe a certain error tolerance and the rank of the solution will adapt at each time step to keep the solution within the level of allowed error. Cross-DEIM builds on ideas developed in \cite{goreinov1997theory} that allow a reconstruction of an $m \times n$ matrix from $r$ of its rows and columns with a pseudo-skeleton approximation. This method is extended to tensors in \cite{oseledets2008tucker} and a Cross$^2$-DEIM algorithm has been developed in \cite{TuckerAA}, however this paper deals only with the two-dimensional case.

The rest of this paper is organized as follows. We give an overview of the OFT in \cref{sec:OFT}, describe the discretization of the Schr\"odinger equation in \cref{sec:low-rank}, and describe its low-rank implementation in \cref{sec:lowrank_impl}. We show experimental results in \cref{sec:experiments}, and our conclusions are in \cref{sec:conclusions}.

\section{Operator Fourier Transform}
\label{sec:OFT}
Following \cite{Cubillos2024OFT}, the OFT defines a function $f$ of an operator $A$ through the Fourier inversion formula. For a scalar argument $\xi$, a function $f$ can be recovered from its Fourier transform $\hat{f}$ as
\begin{equation}
    \label{eq:fourier_inversion}
    f(\xi) = \frac{1}{\sqrt{2\pi}} \int_{-\infty}^{\infty} \hat{f}(\tau)\, e^{i\tau \xi}\, d\tau,
    \qquad
    \hat{f}(\tau) = \frac{1}{\sqrt{2\pi}} \int_{-\infty}^{\infty} f(\xi)\, e^{-i\tau \xi}\, d\xi.
\end{equation}
Replacing the scalar $\xi$ by the operator $A$ and applying the result to a function $g$ gives
\begin{equation}
    \label{eq:oft_full}
    [f(A)g](\vect{x}) = \frac{1}{\sqrt{2\pi}} \int_{-\infty}^{\infty} \hat{f}(\tau)\, e^{i\tau A} g(\vect{x})\, d\tau.
\end{equation}
For a diagonalizable matrix $A = S \Lambda S^{-1}$, as arises after discretization, we have $e^{i\tau A} = S e^{i\tau\Lambda} S^{-1}$, so \cref{eq:oft_full} is simply \cref{eq:fourier_inversion} applied to each eigenvalue of $A$. For all functions $f$ considered in this paper, $\hat{f}(\tau) = 0$ for $\tau < 0$, i.e.\ $\hat{f}$ contains the Heaviside factor $H(\tau)$; see \cref{eq:oft2,eq:f_hat}. The integral in \cref{eq:oft_full} therefore reduces to
\begin{equation}
    \label{eq:oft}
    [f(A)g](\vect{x}) = \frac{1}{\sqrt{2\pi}} \int_{0}^{\infty} \hat{f}(\tau)\, e^{i\tau A} g(\vect{x})\, d\tau,
\end{equation}
which requires $e^{i\tau A}g$ only for nonnegative pseudo-times $\tau$. This is what makes the OFT practical: the function $u(\tau, \vect{x}) \coloneqq e^{i\tau A}g(\vect{x})$ is the solution of the initial value problem
\begin{equation}
    \label{eq:schro_cases}
    \begin{aligned}
        u_\tau(\tau, \vect{x}) &= iAu(\tau, \vect{x}), && \tau > 0, \\
        u(0, \vect{x}) &= g(\vect{x}),
    \end{aligned}
\end{equation}
since differentiating $e^{i\tau A}g$ with respect to $\tau$ gives $iAe^{i\tau A}g$, and $e^{0}g = g$; see \cite{Cubillos2024OFT} for details. Thus $u$ is obtained by evolving \cref{eq:schro_cases} forward in $\tau$ from the initial data $g$, and never backward. In this work \cref{eq:schro_cases} is solved with a low-rank method. Using the OFT framework, the problem of inverting an operator is reduced to solving \cref{eq:schro_cases} and evaluating the integral \cref{eq:oft} numerically. We can expect to approximate the integral quickly when $\hat{f}$ decays quickly in $\tau$.

\subsection{Using the OFT to solve an ODE}
Before we use the OFT to solve \cref{eq:helmholtz}, we revisit a simple example from \cite{Cubillos2024OFT}. Consider the ODE
\begin{equation}
    v(x) - iv_{xx}(x) = g(x).
    \label{eq:oft_ex}
\end{equation}
To use the OFT to solve this equation, we let $A = \partial_{xx}$ and write \cref{eq:oft_ex} as
\begin{equation}
    (I - iA)\, v(x) = g(x),
\end{equation}
so that
\begin{equation}
    v(x) = (I - iA)^{-1} g(x).
\end{equation}
Defining $f(\xi) \coloneqq (1 - i\xi)^{-1}$, so that $f(A) = (I - iA)^{-1}$, we can write the solution as
\begin{equation}\label{eq:oft2}
    v(x) = [f(A)g](x) = \frac{1}{\sqrt{2\pi}} \int_0^\infty \hat{f}(\tau)\, e^{i\tau A} g(x)\, d\tau = \int_{0}^\infty e^{-\tau} e^{i \tau A} g(x) \,d \tau,
\end{equation}
where $\hat{f}(\tau) = \sqrt{2 \pi}\,e^{-\tau}H(\tau)$ is the Fourier transform of $f$ and $H(\tau)$ is the Heaviside function. The representation \cref{eq:oft2} can also be seen directly from the Laplace transform. If $u(\tau, x) = e^{i\tau A}g(x)$ solves \cref{eq:schro_cases}, its Laplace transform $\hat{u}(s, x) = \int_0^\infty e^{-s\tau} u(\tau, x)\, d\tau$ satisfies $s\hat{u} - g = iA\hat{u}$. Evaluating at $s = 1$ gives $(I - iA)\hat{u}(1, x) = g(x)$, so $v = \hat{u}(1, \cdot)$, which is exactly \cref{eq:oft2}. The solution of the initial value problem \cref{eq:schro_cases} can then be inserted into \cref{eq:oft2}, which can in turn be evaluated with a numerical integration scheme.

\subsection{Using the OFT to solve the Helmholtz equation}
\label{sec:oft_helmholtz}

As our objective is to find the solution to the Helmholtz equation, we must first formulate a Helmholtz operator in order to treat it with the OFT. Following \cite{Cubillos2024OFT}, we write \cref{eq:helmholtz} as $Av = g$ with the Helmholtz operator
\begin{equation}
    \label{eq:helmholtz_A}
    A = m(\vect{x}) + \frac{\Delta}{\kappa^2},
\end{equation}
so that $v = A^{-1}g$. Rather than applying the OFT to $f(\xi) = \xi^{-1}$ directly, Cubillos and Jim\'enez take $f(\xi) = \xi^{-1/2}$, so that $f(A)^2 = A^{-1}$. In their words, ``the problem of inverting the Helmholtz operator is reduced to two sequential solves of a pseudo-time paraxial (or Schr\"odinger) equation'' \cite[\S 1]{Cubillos2024OFT}, i.e.\ the solution is obtained by applying $f(A)$ twice, $v = f(A)\big(f(A)g\big)$, and this composition is verified to be the inverse of $A$ in \cite[\S 4]{Cubillos2024OFT}. This choice is convenient because, as we show below, the Fourier transform of $\xi^{-1/2}$ vanishes for $\tau < 0$ and decays like $\tau^{-1/2}$, so only forward solves of \cref{eq:schro_cases} are needed, whereas the Fourier transform of $\xi^{-1}$ is proportional to $\mathrm{sgn}(\tau)$ \cite[\S 4]{Cubillos2024OFT} and has neither property. Throughout, the square root of a complex number is taken to be the principal square root.

Since the Fourier transform of $f(\xi) = \xi^{-1/2}$ is
\begin{equation}
\label{eq:f_hat}
    \hat{f}(\tau) = (1 - i)\,H(\tau)\frac{1}{\sqrt{\tau}},
\end{equation}
we can write the inverse square root Helmholtz operator as
\begin{equation}
    \label{eq:helmholtz_operator}
    f(A) = \frac{1}{\sqrt{2\pi}} \int_0^\infty \hat{f}(\tau)\, e^{i\tau A}\, d\tau = \sqrt{\frac{-i}{\pi}} \int_0^\infty \frac{e^{i \tau}}{\sqrt{\tau}}\, e^{i \tau (A-I)}\, d\tau.
\end{equation}
In the last step we used $(1-i)/\sqrt{2\pi} = \sqrt{-i/\pi}$ and split off the factor $e^{i\tau}$, so that the remaining term $e^{i\tau(A-I)}$ involves the coefficient $m - 1$; this is the form used in \cref{sec:low-rank}. To find the solution to the Helmholtz equation \cref{eq:helmholtz} we apply $f(A)$ to $g$ twice,
\begin{equation}
    \label{eq:helmholtz_oft1}
    \tilde{v} = f(A)g = \sqrt{\frac{-i}{\pi}} \int_0^\infty \frac{e^{i \tau}}{\sqrt{\tau}}\, e^{i \tau (A-I)}g(\vect{x})\, d\tau
\end{equation}
followed by
\begin{equation}
    \label{eq:helmholtz_oft2}
    v = f(A) \tilde{v} = \sqrt{\frac{-i}{\pi}} \int_0^\infty \frac{e^{i \tau}}{\sqrt{\tau}}\, e^{i \tau (A-I)} \tilde{v} (\vect{x})\, d\tau.
\end{equation}

\subsection{Numerical integration of the OFT integral}
\label{sec:num_int_oft}
To numerically integrate the action of $f(A)$ on $g$, we closely follow \cite{Cubillos2024OFT} and use an exponential time-stepping scheme together with the same quadrature weights. We truncate the integral to a finite pseudo-time interval $[0, T]$ and introduce the nodes
\begin{equation}
\label{eq:tau_def}
    \tau_k = a\left(b^k - 1\right), \qquad k = 0, 1, \dots, N,
\end{equation}
so that $\tau_0 = 0$ and the step sizes $\tau_{k+1} - \tau_k = a(b-1)b^k$ grow geometrically by the factor $b > 1$. This allows us to prescribe the step sizes at the start and at the end of the integration, which we denote by $\Delta t_0$ and $\Delta t_f$, respectively. To do so, we let
\begin{equation}
    \label{eq:tau_params}
    R = \frac{\Delta t_f}{\Delta t_0} - 1, \qquad a = \frac{T}{R}, \qquad b = 1 + \frac{R\, \Delta t_0}{T}.
\end{equation}
Here $R + 1$ is the ratio of the final to the initial step size, $b$ is the growth factor between consecutive step sizes, and $a$ sets the overall scale of the nodes. With these choices the first step is $\tau_1 - \tau_0 = a(b-1) = \Delta t_0$. The number of steps $N$ is the smallest integer for which $a(b^N - 1) \geq T$, i.e.\ $N = \lceil \ln(R+1)/\ln b \rceil$, and the last node is set to $\tau_N = T$. The step sizes therefore grow geometrically from $\Delta t_0$ and never exceed $\Delta t_f$.

Following \cite{Cubillos2024OFT}, let $u(\tau, \vect{x}) = e^{i\tau(A-I)}g(\vect{x})$ and $u_k = u(\tau_k, \cdot)$. We replace $u$ by its piecewise linear interpolant on the nodes \cref{eq:tau_def} and integrate the singular kernel $e^{i\tau}/\sqrt{\tau}$ against it exactly. This gives the approximation of \cref{eq:helmholtz_oft1}
\begin{equation}
    \label{eq:oft_quadrature}
    \sqrt{\frac{-i}{\pi}} \int_0^T \frac{e^{i \tau}}{\sqrt{\tau}}\, u(\tau, \vect{x})\, d\tau \approx \sum^N_{k=0}\omega_k u_k,
\end{equation}
and similarly for \cref{eq:helmholtz_oft2}, with the quadrature weights
\begin{equation}
    \label{eq:OFTWeights}
    \omega_k =
    \begin{cases}
        \psi_1\left(\tau_0, \tau_1\right), & k = 0, \\
        \psi_2\left(\tau_{k-1},\tau_k\right) + \psi_1\left(\tau_k,\tau_{k+1}\right), & 1 \leq k \leq N-1, \\
        \psi_2\left(\tau_{N-1},\tau_N\right), & k = N.
    \end{cases}
\end{equation}
Here $\psi_1(\tau_a, \tau_b)$ and $\psi_2(\tau_a, \tau_b)$ are the contributions of the interval $[\tau_a, \tau_b]$ to the weights of its left and right endpoints, respectively,
\begin{equation}
  \label{eq:psi_weights}
  \begin{aligned}
  \psi_1(\tau_a,\tau_b) &= \frac{1+i}{\sqrt{2\pi}(\tau_b-\tau_a)}
              \Big[ \sqrt{\tau_b}\, e^{i \tau_b} - \sqrt{\tau_a}\, e^{i \tau_a}
                    - (1+2i \tau_b)\, D(\tau_a, \tau_b) \Big], \\
  \psi_2(\tau_a,\tau_b) &= \frac{1+i}{\sqrt{2\pi}(\tau_b-\tau_a)}
              \Big[ \sqrt{\tau_a}\, e^{i \tau_a} - \sqrt{\tau_b}\, e^{i \tau_b}
                    + (1+2i \tau_a)\, D(\tau_a, \tau_b) \Big],
  \end{aligned}
\end{equation}
with
\begin{equation}
  D(\tau_a, \tau_b) = C(\sqrt{\tau_b}) - C(\sqrt{\tau_a}) + i\big( S(\sqrt{\tau_b}) - S(\sqrt{\tau_a}) \big),
\end{equation}
where $S$ and $C$ are the Fresnel integrals
\begin{equation}
    \label{eq:fresnel}
    S(\xi) = \int_0^\xi \sin(s^2)\, ds, \qquad C(\xi) = \int_0^\xi \cos(s^2)\, ds.
\end{equation}
We evaluate both through the complex error function,
\begin{equation}
    C(\xi) + iS(\xi) = \frac{\sqrt{\pi}}{2}\, e^{i\pi/4} \operatorname{erf}\!\left(e^{-i\pi/4}\xi\right),
\end{equation}
using the implementation of $\operatorname{erf}$ for complex arguments in the Julia package \texttt{SpecialFunctions.jl}.

\section{Method for finding an approximation to the solution of the Schr\"odinger equation}
\label{sec:low-rank}

To evaluate \cref{eq:helmholtz_oft1} and \cref{eq:helmholtz_oft2} we must first develop a method to evaluate $X(\vect{x}, t) = e^{it(A-I)}g(\vect{x})$. We therefore design a low-rank scheme to evaluate equations of the form
\begin{equation}
    \label{eq:base_eq}
    \begin{aligned}
        X_t &= i\left(\left(m(\vect{x})-1\right)X + \dfrac{\Delta}{\kappa^2}X\right), && \vect{x} \in \Omega, \\[4pt]
        X(\vect{x}, t) &= 0, && \vect{x} \in \partial\Omega, \\[4pt]
        X(\vect{x}, 0) &= g(\vect{x}), && \vect{x} \in \Omega,
    \end{aligned}
\end{equation}
where $\Delta$ is the Laplacian operator. Equation \cref{eq:base_eq} is a time-dependent linear Schr\"odinger equation with the complex potential $m(\vect{x}) - 1$. Its solution at time $t$ is $X(\cdot, t) = e^{it(A-I)}g$, which is the quantity needed at each quadrature node $\tau_k$ in \cref{eq:oft_quadrature}. The coefficient $m(\vect{x})$, the square of the index of refraction, is in general complex-valued, with the real part determining the local wave speed and the imaginary part accounting for absorption in the medium. We therefore write $m(\vect{x}) - 1 = \alpha(\vect{x}) + i\beta(\vect{x})$, where $\alpha(\vect{x}) \in \mathbb{R}$ is the deviation of the real part of $m$ from its free-space value $1$, and $\beta(\vect{x}) \geq 0$ describes the level of absorption. Equation \cref{eq:base_eq} then becomes
\begin{equation}
    \label{eq:base_eq_split}
    X_t = i\left(\alpha + \frac{\Delta}{\kappa^2}\right)X - \beta X.
\end{equation}
We then write $X = V + iW$ with $V$ and $W$ real-valued functions. Substituting into \cref{eq:base_eq_split} and equating real and imaginary parts yields
\begin{equation}
    \label{eq:Vt_Wt}
    \begin{aligned}
        V_t = -\left(\alpha + \dfrac{\Delta}{\kappa^2}\right)W - \beta V, \\[6pt]
        W_t = \left(\alpha + \dfrac{\Delta}{\kappa^2}\right)V - \beta W.
    \end{aligned}
\end{equation}
In this form, $\beta$ acts as a damping term on each component, while $\alpha + \Delta/\kappa^2$ describes the coupling between $V$ and $W$.

\subsection{Discretization on a Cartesian grid}
\label{sec:discretization}

To approximate the solution of \cref{eq:Vt_Wt} numerically we take $\Omega$ to be a square of side length $L$, which after a translation is $\Omega = [0,L]^2$ ($L = 1$ in Example~1 and $L = 2\pi$ in Example~2), and discretize it with a uniform Cartesian grid of $(n+2) \times (n+2)$ points
\begin{equation}
    \label{eq:grid}
    x_i = i h, \qquad y_j = j h, \qquad 0 \le i, j \le n+1, \qquad h = \frac{L}{n+1},
\end{equation}
where $h$ is the spacing between adjacent grid points in both directions. The boundary $\partial\Omega$ corresponds to the indices $i \in \{0, n+1\}$ and $j \in \{0, n+1\}$, leaving the $n \times n$ interior points $(x_i, y_j)$, $1 \le i, j \le n$, as the unknowns of the discrete problem.

On this grid we introduce the real-valued grid functions
\begin{equation}
    V_{ij}(t) \approx V(x_i, y_j, t), \qquad W_{ij}(t) \approx W(x_i, y_j, t),
\end{equation}
and assemble them into $n \times n$ matrices $V(t), W(t) \in \mathbb{R}^{n \times n}$ with entries $V_{ij}$, $W_{ij}$. The homogeneous Dirichlet condition $X = 0$ on $\partial \Omega$ enforces $V_{0,j} = V_{n+1,j} = V_{i,0} = V_{i,n+1} = 0$ and likewise for $W$, so the boundary values do not appear as unknowns.

Following the standard centered-difference approximation \cite[\S 6.3]{demmel1997book}, the one-dimensional second derivative is approximated at an interior point $x_i$ by
\begin{equation}
    \label{eq:1d_second_deriv}
    -\frac{d^2 V}{dx^2}\bigg|_{x = x_i} \;\approx\; \frac{2V_i - V_{i-1} - V_{i+1}}{h^2}.
\end{equation}
Collecting these approximations for $i = 1, \dots, n$ and using the Dirichlet boundary values $V_0 = V_{n+1} = 0$ gives the symmetric positive-definite tridiagonal matrix
\begin{equation}
    \label{eq:Tn_def}
    T_n = \begin{bmatrix}
        2  & -1 &        &        &    \\
        -1 &  2 & -1     &        &    \\
           & -1 & \ddots & \ddots &    \\
           &    & \ddots & \ddots & -1 \\
           &    &        & -1     &  2
    \end{bmatrix} \in \mathbb{R}^{n \times n},
\end{equation}
so that $-d^2/dx^2 \approx h^{-2} T_n$ in the discrete setting. In two dimensions, applying \cref{eq:1d_second_deriv} along each coordinate direction yields the standard second-order finite-difference approximation of the Laplacian \cite[\S 6.3.2]{demmel1997book}
\begin{equation}
    \label{eq:five_point}
    -\Delta V|_{(x_i, y_j)} \;\approx\; \frac{1}{h^2}\bigl( T_n V + V T_n \bigr)_{ij},
\end{equation}
so the discrete Laplacian acts on the matrix of unknowns directly. The matrix $T_n$ has the eigendecomposition
\begin{equation}
    \label{eq:Tn_eig}
    T_n = Z \Lambda Z^T, \qquad \Lambda = \mathrm{diag}(\lambda_1, \dots, \lambda_n), \qquad \lambda_\ell = 2\left(1 - \cos\frac{\ell\pi}{n+1}\right),
\end{equation}
with orthonormal eigenvectors $Z_{i\ell} = \sqrt{2/(n+1)}\,\sin\bigl(i\ell\pi/(n+1)\bigr)$ \cite[Lemma 6.1]{demmel1997book}, so $Z$ and $Z^T$ can be applied in $O(n \log n)$ operations via the discrete sine transform (DST).

We discretize \cref{eq:Vt_Wt} in time using the midpoint rule. From here on we assume that the coefficient $m$, and therefore $\alpha$ and $\beta$, is constant in $\Omega$, as it is in both examples of \cref{sec:experiments}. This is what allows the transform to the sine basis below: for a constant $\alpha$ we have $Z^T(\alpha \bar{W})Z = \alpha\, Z^T\bar{W}Z$, whereas for a variable coefficient the entry-wise product with $\alpha(x_i, y_j)$ is not diagonalized by $Z$ and the update \cref{eq:step4} below no longer decouples. With $V^k, W^k$ the grid-function matrices at time $t^k$ and $\Delta t_k = t^{k+1} - t^k$ the step size, which varies with $k$ for the exponential time stepping of \cref{sec:num_int_oft}, define
\begin{equation}
    \label{eq:Vbar_Wbar}
    \bar{V} = \tfrac{1}{2}\!\left(V^{k+1} + V^k\right), \qquad
    \bar{W} = \tfrac{1}{2}\!\left(W^{k+1} + W^k\right),
\end{equation}
so that, using \cref{eq:five_point} to discretize $\Delta$,
\begin{equation}
    \label{eq:step1}
    \begin{aligned}
        \frac{V^{k+1} - V^k}{\Delta t_k} &= -\left(\alpha \bar{W} - \frac{1}{h^2 \kappa^2} \left(T_n \bar{W} + \bar{W} T_n\right)\right) - \beta \bar{V}, \\[1ex]
        \frac{W^{k+1} - W^k}{\Delta t_k} &= \left(\alpha \bar{V} - \frac{1}{h^2 \kappa^2} \left(T_n \bar{V} + \bar{V} T_n\right)\right) - \beta \bar{W}.
    \end{aligned}
\end{equation}
Using the eigendecomposition \cref{eq:Tn_eig} we apply $Z$ and $Z^T$ via the DST and rewrite
\begin{equation}
    \label{eq:step2}
    \begin{aligned}
        V^{k+1} - V^k &= \Delta t_k \left(-\alpha \bar{W} + \frac{1}{h^2\kappa^2} \left(Z \Lambda Z^T \bar{W} + \bar{W} Z \Lambda Z^T\right) - \beta \bar{V}\right), \\
        W^{k+1} - W^k &= \Delta t_k \left(\alpha \bar{V} - \frac{1}{h^2\kappa^2} \left(Z \Lambda Z^T \bar{V} + \bar{V} Z \Lambda Z^T\right) - \beta \bar{W}\right).
    \end{aligned}
\end{equation}
With the identifications
\begin{equation}
    \label{eq:identification1}
    V^{k+1} - V^k = 2\bar{V} - 2V^k, \qquad
    W^{k+1} - W^k = 2\bar{W} - 2W^k,
\end{equation}
and defining the transformed quantities
\begin{equation}
    \label{eq:hat_definition}
    \hat{V} = Z^T \bar{V} Z, \quad \hat{V}^k = Z^T V^k Z, \quad
    \hat{W} = Z^T \bar{W} Z, \quad \hat{W}^k = Z^T W^k Z,
\end{equation}
the system becomes
\begin{equation}
    \label{eq:step3}
    \begin{aligned}
        \hat{V} - \hat{V}^k &= \frac{\Delta t_k}{2}\left(-\alpha \hat{W} + \frac{1}{h^2\kappa^2}\left(\Lambda \hat{W} + \hat{W} \Lambda\right) - \beta \hat{V}\right), \\
        \hat{W} - \hat{W}^k &= \frac{\Delta t_k}{2}\left(\alpha \hat{V} - \frac{1}{h^2\kappa^2}\left(\Lambda \hat{V} + \hat{V} \Lambda\right) - \beta \hat{W}\right).
    \end{aligned}
\end{equation}
Since $\Lambda$ is diagonal, $(\Lambda \hat{W} + \hat{W}\Lambda)_{ij} = (\lambda_i + \lambda_j)\hat{W}_{ij}$, and similarly for $\hat{V}$. The $(i,j)$ entry of \cref{eq:step3} therefore involves only $\hat{V}_{ij}$ and $\hat{W}_{ij}$, so the system decouples into $n^2$ independent $2 \times 2$ systems, one for each pair of indices $i, j \in \{1, \dots, n\}$:
\begin{equation}
    \label{eq:step4}
    \begin{bmatrix}
        1 + \tfrac{\beta \Delta t_k}{2} & \tfrac{\Delta t_k}{2}\!\left(\alpha - \tfrac{1}{h^2\kappa^2}\!\left(\lambda_i + \lambda_j\right)\right) \\[6pt]
        -\tfrac{\Delta t_k}{2}\!\left(\alpha - \tfrac{1}{h^2\kappa^2}\!\left(\lambda_i + \lambda_j\right)\right) & 1 + \tfrac{\beta \Delta t_k}{2}
    \end{bmatrix}\!\begin{bmatrix}
        \hat{V}_{ij} \\ \hat{W}_{ij}
    \end{bmatrix} = \begin{bmatrix}
        \hat{V}^k_{ij} \\ \hat{W}^k_{ij}
    \end{bmatrix}.
\end{equation}
Within each system, $\hat{V}_{ij}$ and $\hat{W}_{ij}$ remain coupled; what decouples are the different grid points.
Defining
\begin{equation}
    \label{eq:c_and_d}
    c = 1 + \frac{\beta \Delta t_k}{2}, \qquad
    d = \frac{\Delta t_k}{2}\left(\alpha - \frac{1}{h^2\kappa^2}\left(\lambda_i + \lambda_j\right)\right),
\end{equation}
we may rewrite \cref{eq:step4} as
\begin{equation}
\label{eq:c_d_intermediate}
    \begin{bmatrix}
        c & d \\ -d & c
    \end{bmatrix}\begin{bmatrix}
        \hat{V}_{ij} \\ \hat{W}_{ij}
    \end{bmatrix} = \begin{bmatrix}
        \hat{V}^k_{ij} \\ \hat{W}^k_{ij}
    \end{bmatrix},
\end{equation}
which by Cramer's rule gives the midpoint values
\begin{equation}
    \hat{V}_{ij} = \frac{c \hat{V}^k_{ij} - d \hat{W}^k_{ij}}{c^2 + d^2}, \qquad
    \hat{W}_{ij} = \frac{d \hat{V}^k_{ij} + c \hat{W}^k_{ij}}{c^2 + d^2}.
\end{equation}
By linearity of the discrete sine transform, $\hat{V}$ and $\hat{W}$ are averages of $\hat{V}^{k+1}$ and $\hat{V}^k$ and of $\hat{W}^{k+1}$ and $\hat{W}^k$, that is $\hat{V}^{k+1}_{ij} = 2\hat{V}_{ij} - \hat{V}^k_{ij}$ (and similarly for $\hat{W}$). Since we are interested in computing $\hat{V}^{k+1}$ and $\hat{W}^{k+1}$ we now express the full system in terms of non-averaged quantities
\begin{equation}
    \label{eq:two_sided_system}
    \begin{bmatrix}
        c & d \\ -d & c
    \end{bmatrix}\begin{bmatrix}
        \hat{V}^{k+1}_{ij} \\ \hat{W}^{k+1}_{ij}
    \end{bmatrix} = \begin{bmatrix}
        2 - c & -d \\ d & 2 - c
    \end{bmatrix}\begin{bmatrix}
        \hat{V}^k_{ij} \\ \hat{W}^k_{ij}
    \end{bmatrix}.
\end{equation}
Inverting the left-hand side, the direct update rule is
\begin{equation}
    \label{eq:update_rule}
    \begin{bmatrix}
        \hat{V}_{ij}^{k+1} \\ \hat{W}_{ij}^{k+1}
    \end{bmatrix} = \frac{1}{c^2 + d^2}\begin{bmatrix}
        c\left(2 - c\right) - d^2 & -2d \\
        2d & c\left(2 - c\right) - d^2
    \end{bmatrix}\begin{bmatrix}
        \hat{V}^k_{ij} \\ \hat{W}^k_{ij}
    \end{bmatrix}.
\end{equation}
The update matrix in \cref{eq:update_rule} has the form $\frac{1}{c^2 + d^2}\left(pI + qJ\right)$ where $J = \begin{bmatrix} 0 & -1 \\ 1 & 0 \end{bmatrix}$ and
\begin{equation}
    p = c\left(2 - c\right) - d^2, \qquad q = 2d.
\end{equation}
Since any matrix of the form $pI + qJ$ can be decomposed as $\sqrt{p^2 + q^2}\begin{bmatrix} \cos\theta & -\sin\theta \\ \sin\theta & \cos\theta \end{bmatrix}$, we can write
\begin{equation}
    \label{eq:theta_update}
    \begin{bmatrix}
        \hat{V}_{ij}^{k+1} \\ \hat{W}_{ij}^{k+1}
    \end{bmatrix} = \phi\begin{bmatrix}
        \cos\theta & -\sin\theta \\ \sin\theta & \cos\theta
    \end{bmatrix}\begin{bmatrix}
        \hat{V}^k_{ij} \\ \hat{W}^k_{ij}
    \end{bmatrix}
\end{equation}
where
\begin{equation}
    \label{eq:theta_rho_def}
    \theta = \operatorname{atan2}(q, p), \qquad \phi = \frac{\sqrt{p^2 + q^2}}{c^2 + d^2}.
\end{equation}
Here $\operatorname{atan2}(q, p)$ is the angle of the point $(p, q)$, i.e.\ $\cos\theta = p/\sqrt{p^2 + q^2}$ and $\sin\theta = q/\sqrt{p^2 + q^2}$; the two-argument form is needed because $p$ can be negative. The prefactor $\phi$ controls the amplitude change at each step. When $\beta = 0$ we have $c = 1$, so $p = 1 - d^2$, $q = 2d$, and $p^2 + q^2 = \left(1 + d^2\right)^2 = \left(c^2 + d^2\right)^2$, giving $\phi = 1$. In this case the scheme reduces to a pure rotation
\begin{equation}
    \begin{bmatrix}
        \hat{V}_{ij}^{k+1} \\ \hat{W}_{ij}^{k+1}
    \end{bmatrix} = \begin{bmatrix}
        \cos\theta & -\sin\theta \\ \sin\theta & \cos\theta
    \end{bmatrix}\begin{bmatrix}
        \hat{V}^k_{ij} \\ \hat{W}^k_{ij}
    \end{bmatrix}
\end{equation}
with $\theta = 2\arctan(d)$, and the scheme is norm-preserving. When $\beta > 0$, we have $\phi < 1$, and the update introduces damping at each step, consistent with the dissipative terms $-\beta V$, $-\beta W$ in \cref{eq:Vt_Wt}. Once $\hat{V}^{k+1}$ and $\hat{W}^{k+1}$ have been computed in the transformed basis via \cref{eq:theta_update}, the solution in the original variables is recovered by inverting the transform in \cref{eq:hat_definition},
\begin{equation}
    \label{eq:inverse_dst}
    V^{k+1} = Z \hat{V}^{k+1} Z^T, \qquad W^{k+1} = Z \hat{W}^{k+1} Z^T,
\end{equation}
which, like the forward transform, can be applied in $O(n^2 \log n)$ operations via the discrete sine transform.

\section{Low-rank implementation of the OFT}
\label{sec:lowrank_impl}
In this section we describe how the time step of \cref{sec:discretization} and the quadrature of \cref{sec:num_int_oft} are carried out when $V$ and $W$ are stored in low-rank form. \Cref{sec:truncsum,sec:cross-deim} review the two low-rank building blocks we use: the rounded sum of low-rank matrices and the Cross-DEIM algorithm. \Cref{sec:lr_timestep} then describes the low-rank time step, and \cref{sec:lr_oft} the full low-rank OFT solver and its cost.

\subsection{Low-rank arithmetic: rounded sum of matrices}
\label{sec:truncsum}
During each step of a Schr\"odinger solve, and when accumulating the OFT quadrature sum, we must compute sums of low-rank matrices. The sum of low-rank matrices is in general of higher rank than the summands, so it must be recompressed after it is computed. Let $C_1, \dots, C_s \in \mathbb{R}^{m \times n}$ be given in SVD form, $C_i = U_i \Sigma_i V_i^T$ with $U_i \in \mathbb{R}^{m \times r_i}$, $\Sigma_i \in \mathbb{R}^{r_i \times r_i}$ and $V_i \in \mathbb{R}^{n \times r_i}$, and let $r_{\rm tot} = \sum_{i=1}^s r_i$ be the total rank. We want to compute a truncated SVD of
\begin{equation}
    X = \sum_{i=1}^s C_i = \sum_{i=1}^s U_i \Sigma_i V_i^T.
\end{equation}
We first write the sum as a single product,
\begin{equation}
    \sum_{i=1}^s U_i \Sigma_i V_i^T = \begin{bmatrix}
    U_1 & \cdots & U_s
\end{bmatrix}\begin{bmatrix}
    \Sigma_1 & & \\ & \ddots & \\ & & \Sigma_s
\end{bmatrix}\begin{bmatrix}
    V_1^T \\ \vdots \\ V^T_s
\end{bmatrix} = U_{big}\Sigma_{big}V^T_{big},
\end{equation}
where $U_{big} \in \mathbb{R}^{m \times r_{\rm tot}}$, $\Sigma_{big} \in \mathbb{R}^{r_{\rm tot} \times r_{\rm tot}}$ and $V_{big} \in \mathbb{R}^{n \times r_{\rm tot}}$. In general neither $U_{big}$ nor $V_{big}$ has orthonormal columns, so this is not an SVD of $X$. To obtain orthonormal bases we compute the column pivoted QR decompositions
\begin{equation}
    U_{big} = Q_U R_U \Pi_U^T, \qquad V_{big} = Q_V R_V \Pi_V^T,
\end{equation}
where $Q_U \in \mathbb{R}^{m \times r_{\rm tot}}$ and $Q_V \in \mathbb{R}^{n \times r_{\rm tot}}$ have orthonormal columns, $R_U, R_V \in \mathbb{R}^{r_{\rm tot} \times r_{\rm tot}}$ are upper triangular, and $\Pi_U, \Pi_V$ are permutation matrices. We then form the small matrix
\begin{equation}
    L = R_U \Pi_U^T \Sigma_{big}\Pi_V R_V^T \in \mathbb{R}^{r_{\rm tot} \times r_{\rm tot}},
\end{equation}
so that $X = Q_U L Q_V^T$. Computing the truncated SVD $L \approx U_L \Sigma_L V_L^T$, in which singular values below a prescribed tolerance are discarded, gives
\begin{equation}
    X \approx U_X\Sigma_XV^T_X, \qquad U_X = Q_U U_L, \quad \Sigma_X = \Sigma_L, \quad V_X = Q_V V_L.
\end{equation}
We refer to this procedure as \texttt{truncsum}; it is summarized in \cref{alg:truncsum}.

All operations above act on matrices with at most $r_{\rm tot}$ columns. The two QR decompositions cost $O((m+n)r_{\rm tot}^2)$, forming and decomposing $L$ costs $O(r_{\rm tot}^3)$, and forming $U_X$ and $V_X$ costs $O((m+n)r_{\rm tot}^2)$, for a total of $O((m+n)r_{\rm tot}^2 + r_{\rm tot}^3)$. In contrast, assembling $X$ as a dense matrix costs $O(mnr_{\rm tot})$, and computing its SVD costs $O(mn\min(m,n))$. A column pivoted QR decomposition of the dense matrix $X$ is no cheaper: it costs $O(mn\min(m,n))$, or $O(mnr)$ if it is stopped at the truncation rank $r$, in addition to the cost of forming $X$. Since $r_{\rm tot} \ll \min(m,n)$ in our setting, the rounded sum is much cheaper, and the dense matrix $X$ is never formed.

\begin{algorithm}[]
\caption{Rounding of sum of low-rank matrices, $U_X\Sigma_X V_X^T = \mathcal{T}^{\rm round}_{\epsilon,r_{\rm max}}(\sum_{i=1}^s U_i\Sigma_iV^T_i)$ \label{alg:truncsum}}

\begin{algorithmic}[]
\STATE{{\bf Input:} Low-rank matrices in the form $U_i\Sigma_iV^T_i$, $i =1,\dots,s$, tolerance $\epsilon$, max rank $r_{\rm max}$}
\STATE{{\bf Output:} $U_X,\Sigma_X,V_X$}
\STATE{Let $U_{big}=[U_1,\dots ,U_s]$, $\Sigma_{big}={\sf diag}(\Sigma_1,\dots,\Sigma_s)$, $V_{big}=[V_1,\dots,V_s]$}
\STATE{Perform column pivoted QR: $[Q_U,R_U,\Pi_U] = {\sf qr}(U_{big})$, $[Q_V,R_V,\Pi_V] = {\sf qr}(V_{big})$}
\STATE{Compute the truncated SVD of the small matrix with tolerance $\epsilon$ and max rank $r_{\rm max}$: \[\mathcal{T}_{\epsilon,r_{\rm max}}(R_U\Pi_U^T\Sigma_{big}\Pi_V R_V^T) = U_L\Sigma_L V_L^T\]}
\STATE{Form $U_X \leftarrow Q_U U_L$, $\Sigma_X \leftarrow \Sigma_L$, $V_X \leftarrow Q_V V_L$}
\RETURN $[U_X, \Sigma_X, V_X]$
\end{algorithmic}
\end{algorithm}

\subsection{Cross-DEIM}
\label{sec:cross-deim}
In this section we describe how to compute a low-rank SVD approximation of a matrix $G \in \mathbb{R}^{m \times n}$ with the Cross-DEIM algorithm introduced in \cite{appelo2025lraa}. We assume only that individual entries $G(i,j)$ can be evaluated, for example through a function handle; the full matrix $G$ is never formed. Cross-DEIM samples only a few rows and columns of $G$, so its cost scales sub-linearly in the number of entries $mn$, whereas a truncated SVD requires all of $G$. In outline, Cross-DEIM proceeds as follows:
\begin{enumerate}
    \item select a set $\mathcal{I}$ of row indices and a set $\mathcal{J}$ of column indices;
    \item form the small matrix $G(\mathcal{I}, \mathcal{J})$ and $M = G(\mathcal{I}, \mathcal{J})^{\dagger}$;
    \item compute the SVD $U S V^T = G(:, \mathcal{J})\, M\, G(\mathcal{I}, :)$, which requires only the columns $G(:, \mathcal{J})$ and the rows $G(\mathcal{I}, :)$;
    \item check the accuracy of the approximation; if it is not sufficient, use the new singular vectors $U$ and $V$ to select better index sets and return to step 2.
\end{enumerate}
We now describe these steps in turn.

\paragraph{Cross approximation (steps 2 and 3)}
Given the index sets, cross approximation approximates $G$ by
\begin{equation}
    \label{eq:cross_approx}
    G \approx G(:, \mathcal{J})\, M\, G(\mathcal{I}, :), \qquad M = G(\mathcal{I}, \mathcal{J})^{\dagger},
\end{equation}
where $G(:, \mathcal{J}) \in \mathbb{R}^{m \times |\mathcal{J}|}$ and $G(\mathcal{I}, :) \in \mathbb{R}^{|\mathcal{I}| \times n}$. This approximation reproduces $G$ exactly in the selected rows and columns \cite{goreinov1997theory}. It is similar to the CUR decomposition, but there the middle matrix is $M = G(:,\mathcal{J})^\dagger\, G\, G(\mathcal{I},:)^\dagger$, which requires access to all of $G$ \cite{mahoney2009cur}. Since $G(\mathcal{I},\mathcal{J})$ can be ill-conditioned, the pseudo-inverse is not formed explicitly but is computed stably using the methods of \cite{donello2023oblique, appelo2025lraa}. Given a function handle to $G(i,j)$ and the two index sets, the subroutine \texttt{scross} in \cref{algo:scross} returns the approximation \cref{eq:cross_approx} in SVD form, $USV^T \approx G$. When $|\mathcal{I}|, |\mathcal{J}| = O(r)$, its cost scales as $O((m+n)r^2)$.

\paragraph{Index selection (step 1)}
The quality of \cref{eq:cross_approx} depends on the choice of $\mathcal{I}$ and $\mathcal{J}$, and choosing them optimally would require the full matrix $G$. Instead, given matrices $U \in \mathbb{R}^{m \times r}$ and $V \in \mathbb{R}^{n \times r}$ with orthonormal columns that approximate the leading singular vectors of $G$, we select the indices greedily. Possible choices include DEIM \cite{DEIMCUR}, QDEIM \cite{drmac2016new}, and leverage scores \cite{mahoney2009cur}; here we use QDEIM. QDEIM computes the column pivoted QR decomposition
\begin{equation}
    \label{eq:qdeim}
    U^T P = Q R, \qquad \mathcal{I} = \{p_1, \dots, p_r\},
\end{equation}
where $p_1, \dots, p_r$ are the first $r$ pivot indices, i.e.\ the rows of $U$ selected first by the pivoting. This choice keeps the submatrix $U(\mathcal{I}, :)$ well-conditioned \cite{drmac2016new}. The column indices $\mathcal{J}$ are obtained in the same way from $V$. The procedure is given in \cref{algo:qdeim}.

\paragraph{Iteration and stopping criterion (step 4)}
Index selection requires the singular vectors of $G$, which are exactly what we want to compute. Cross-DEIM resolves this by iteration. It starts from an initial guess $U_0, V_0$ (in our solver, the singular vectors from the previous time step) and empty index sets $\mathcal{I}_0 = \mathcal{J}_0 = \emptyset$. At iteration $k$ the index sets are enriched with the indices selected by QDEIM from the current singular vectors,
\begin{equation}
    \mathcal{I}_k = \mathcal{I}_{k-1} \cup \texttt{QDEIM}(U_{k-1}), \qquad \mathcal{J}_k = \mathcal{J}_{k-1} \cup \texttt{QDEIM}(V_{k-1}),
\end{equation}
with a randomly chosen index added if no new index is found, and the stabilized cross approximation gives a new approximation $U_k S_k V_k^T = \texttt{scross}(G, \mathcal{I}_k, \mathcal{J}_k)$. Rows and columns that are numerically linearly dependent, which are detected from the diagonals of the pivoted QR factors in \texttt{scross}, are removed from the index sets. The iteration stops when both the change between consecutive approximations, $\| U_k S_k V_k^T - U_{k-1} S_{k-1} V_{k-1}^T \|$, and an error indicator from \cite{donello2023oblique}, based on the smallest retained singular value and on the conditioning of $U_k(\mathcal{I}_k,:)$ and $V_k(\mathcal{J}_k,:)$, are below the tolerance $\epsilon$. Finally, the approximation is truncated to the smallest rank $r^\ast$ for which the discarded singular values satisfy $\sum_{l > r^\ast} S_l^2 < \epsilon^2$. The full algorithm is given in \cref{algo:crossDEIM}. Assuming that the number of iterations is $O(1)$ and that the index sets have $O(r)$ elements, the cost of Cross-DEIM is $O((m+n)r^2)$ plus the cost of evaluating the $O((m+n)r)$ sampled entries of $G$.

\begin{algorithm}[]
\caption{{\tt [$\mathcal{I}$] = QDEIM(U)} \\QDEIM index selection}
 \begin{algorithmic}[1]
 \STATE {\bf Input:}  Matrix $U \in \mathbb{R}^{m \times r}$ with orthonormal columns
\STATE {\bf Output:} Index set $\mathcal{I}$ of size $r$
\STATE  $[\sim,\sim, p]=\textrm{qr}(U^T,\textrm{'vector'})$ \COMMENT{Perform column pivoted QR on $U^T$.}
\STATE $\mathcal{I}=p(1:r)$
\RETURN $\mathcal{I}$.
\end{algorithmic}
\label{algo:qdeim}
\end{algorithm}

\begin{algorithm}[]
\caption{{\tt [U,S,V,rC,rR] = scross(G,I,J)} \\Stabilized cross  approximation of $G$}
 \begin{algorithmic}[1]
 \STATE {\bf Input:}  $G$, and two index sets $\mathcal{I}$, $\mathcal{J}$ with $k$ and $l$ elements, respectively.
\STATE {\bf Output:} Approximate SVD of $G$, $U\in \mathbb{R}^{m \times r},$ $S \in \mathbb{R}^{r \times r}$, $V\in \mathbb{R}^{n \times r}$, and two vectors $r_{\rm R}, r_{\rm C}$ to test for linear dependence.
\STATE $C = G(:,\mathcal{J}) \in \mathbb{R}^{m \times l}, \ \ R = G(\mathcal{I},:) \in \mathbb{R}^{k \times n}$
\STATE $CP_{\rm C} = Q_{\rm C}R_{\rm C}$, \ \ $R^TP_{\rm R} = Q_{\rm R}R_{\rm R}$ \COMMENT{Perform column pivoted QR.}
\IF{$l \le k$}
\STATE Solve $Q_{\rm C}(\mathcal{I},:) W = R$ \COMMENT{Solved using {\tt \textbackslash} if $Q_{\rm C}(\mathcal{I},:)$ is well conditioned,}
\STATE \COMMENT{else solved by truncated SVD pseudoinverse.}
\STATE $W = \hat{U}SV^T$  \COMMENT{Perform truncated SVD.}
\STATE $U = Q_{\rm C}\hat{U}$
\ELSE
\STATE Solve $Q_{\rm R}(\mathcal{J},:) W = C^T$ \COMMENT{Solved using {\tt \textbackslash} if $Q_{\rm R}(\mathcal{J},:)$ is well conditioned,}
\STATE  \COMMENT{else solved by truncated SVD pseudoinverse.}
\STATE $W^T = US\hat{V}^T$  \COMMENT{Perform truncated SVD.}
\STATE $V = Q_{\rm R}\hat{V}$
\ENDIF
\RETURN $USV^T \approx G$ and $r_{\rm R} = {\sf diag}(P_R^T R_R P_R), r_{\rm C} = {\sf diag} (P_C^T R_C P_C)$.
\end{algorithmic}
\label{algo:scross}
\end{algorithm}

\subsection{Low-rank time stepping of the Schr\"{o}dinger equation}
\label{sec:lr_timestep}
We now describe the low-rank version of the time step of \cref{sec:discretization}. At time $t^k$ both $V^k$ and $W^k$ are stored as truncated SVDs,
\begin{equation}
    \label{eq:lr_factors}
    V^k = U_V \Sigma_V V_V^T, \qquad W^k = U_W \Sigma_W V_W^T,
\end{equation}
with $U_V, V_V \in \mathbb{R}^{n \times r_V}$, $\Sigma_V \in \mathbb{R}^{r_V \times r_V}$, and similarly for $W$ with rank $r_W$. The dense $n \times n$ matrices are never formed. A time step consists of three stages.

\paragraph{1. Transform to the sine basis}
By \cref{eq:hat_definition}, the transformed matrices are
\begin{equation}
    \label{eq:lr_dst}
    \hat{V}^k = Z^T V^k Z = \hat{U}_V \Sigma_V \hat{V}_V^T, \qquad \hat{U}_V = Z^T U_V, \quad \hat{V}_V = Z^T V_V,
\end{equation}
and similarly $\hat{W}^k = \hat{U}_W \Sigma_W \hat{V}_W^T$. The transform therefore acts only on the factors: computing $\hat{U}_V$ and $\hat{V}_V$ requires $2r_V$ discrete sine transforms of length $n$, at a cost of $O(r_V\, n \log n)$.

\paragraph{2. Entry-wise update with Cross-DEIM}
By \cref{eq:theta_update}, the transformed matrices at the next time step have entries
\begin{equation}
    \label{eq:lr_entry_functions}
    \begin{aligned}
        \hat{V}^{k+1}_{ij} = G_V(i,j) &\coloneqq \phi_{ij}\left(\cos\theta_{ij}\, \hat{V}^k_{ij} - \sin\theta_{ij}\, \hat{W}^k_{ij}\right), \\
        \hat{W}^{k+1}_{ij} = G_W(i,j) &\coloneqq \phi_{ij}\left(\sin\theta_{ij}\, \hat{V}^k_{ij} + \cos\theta_{ij}\, \hat{W}^k_{ij}\right),
    \end{aligned}
\end{equation}
where $\theta_{ij}$ and $\phi_{ij}$ are given by \cref{eq:theta_rho_def} and depend on $(i,j)$ through $\lambda_i + \lambda_j$. A single entry of $\hat{V}^k$ is evaluated from the factors in \cref{eq:lr_dst} as
\begin{equation}
    \label{eq:lr_entry}
    \hat{V}^k_{ij} = \hat{U}_V(i,:)\, \Sigma_V\, \hat{V}_V(j,:)^T,
\end{equation}
at a cost of $O(r_V)$, and similarly for $\hat{W}^k_{ij}$. Each entry of $G_V$ and $G_W$ can therefore be evaluated in $O(r_V + r_W)$ operations. Because $\theta_{ij}$ and $\phi_{ij}$ vary with $(i,j)$, the update \cref{eq:lr_entry_functions} cannot be written as a product of the factors, so we use Cross-DEIM (\cref{sec:cross-deim}) to compute low-rank approximations
\begin{equation}
    \label{eq:lr_crossdeim}
    \begin{aligned}
        \hat{V}^{k+1} &\approx \hat{U}^{k+1}_V \Sigma^{k+1}_V (\hat{V}^{k+1}_V)^T = \texttt{Cross-DEIM}(G_V), \\
        \hat{W}^{k+1} &\approx \hat{U}^{k+1}_W \Sigma^{k+1}_W (\hat{V}^{k+1}_W)^T = \texttt{Cross-DEIM}(G_W),
    \end{aligned}
\end{equation}
with tolerance $\Delta t_k^3$, maximum rank $n$, and the factors of $\hat{V}^k$ as the initial guess; at most 20 iterations are allowed ($\texttt{maxiter} = 20$ in \cref{algo:crossDEIM}), and the parameter $\aleph_{\rm max}$ is omitted from the notation. Both $G_V$ and $G_W$ are evaluated from the factors at time $t^k$, so the new factors are stored separately until both approximations have been computed. Cross-DEIM only evaluates $G_V$ and $G_W$ at the rows and columns it selects.

\paragraph{3. Transform back}
Since $Z$ is orthogonal, inverting \cref{eq:hat_definition} gives
\begin{equation}
    \label{eq:lr_inverse_dst}
    V^{k+1} = Z \hat{V}^{k+1} Z^T = \left(Z \hat{U}^{k+1}_V\right) \Sigma^{k+1}_V \left(Z \hat{V}^{k+1}_V\right)^T,
\end{equation}
and similarly for $W^{k+1}$, so again only the factors are transformed.

The time step is summarized in \cref{alg:timestep}: line 3 performs stage 1, lines 4--10 define $G_W$ and $G_V$ as function handles that are evaluated only at the indices requested by Cross-DEIM, lines 11--12 perform stage 2, and lines 13--14 perform stage 3.

\begin{algorithm}[]
\caption{{\tt [W-factors, V-factors] = LR-Timestep(W-factors, V-factors, $\{\lambda_j\}$, $\Delta t_k$, $\alpha$, $\beta$, $\kappa$)}\\
Low-rank time step for the Schr\"odinger equation via rotation in eigenspace
\label{alg:timestep}}
\begin{algorithmic}[1]
\STATE {\bf Input:} Low-rank factors $U_W\Sigma_W V_W^T$, $U_V\Sigma_V V_V^T$; eigenvalues $\{\lambda_j\}$; step size $\Delta t_k$; parameters $\alpha, \beta, \kappa$
\STATE {\bf Output:} Updated low-rank factors for $W$ and $V$
\STATE $\hat{U}_W \leftarrow \texttt{DST}(U_W), \quad \hat{V}_W \leftarrow \texttt{DST}(V_W), \quad \hat{U}_V \leftarrow \texttt{DST}(U_V), \quad \hat{V}_V \leftarrow \texttt{DST}(V_V)$ \COMMENT{Transform factors via DST.}
\STATE Define entry-wise functions $G_W(i,j)$ and $G_V(i,j)$:
\STATE \quad $w \leftarrow \big(\hat{U}_W(i,:) \Sigma_W\big) \cdot \hat{V}_W(j,:)^T$
\STATE \quad $v \leftarrow \big(\hat{U}_V(i,:) \Sigma_V\big) \cdot \hat{V}_V(j,:)^T$
\STATE \quad $c \leftarrow 1 + \beta \Delta t_k / 2, \quad d \leftarrow \tfrac{\Delta t_k}{2}\!\left(\alpha - \tfrac{\lambda_i + \lambda_j}{h^2\kappa^2}\right)$
\STATE \quad $p \leftarrow c(2 - c) - d^2, \quad q \leftarrow 2d, \quad \theta \leftarrow \operatorname{atan2}(q, p), \quad \phi \leftarrow \sqrt{p^2 + q^2}\,/\,(c^2 + d^2)$
\STATE \quad $G_W(i,j) \leftarrow \phi\!\left(\sin\theta\;v + \cos\theta\;w\right)$
\STATE \quad $G_V(i,j) \leftarrow \phi\!\left(\cos\theta\;v - \sin\theta\;w\right)$
\STATE $[\hat{U}_W^{+}, \Sigma_W^{+}, \hat{V}_W^{+}] \leftarrow \texttt{Cross-DEIM}(G_W, \hat{U}_V, \hat{V}_V, \epsilon = \Delta t_k^3, r_{\rm max} = n)$ \COMMENT{$G_W$ and $G_V$ use the factors from time $t^k$.}
\STATE $[\hat{U}_V^{+}, \Sigma_V^{+}, \hat{V}_V^{+}] \leftarrow \texttt{Cross-DEIM}(G_V, \hat{U}_V, \hat{V}_V, \epsilon = \Delta t_k^3, r_{\rm max} = n)$
\STATE $U_W \leftarrow \texttt{DST}(\hat{U}_W^{+}), \quad V_W \leftarrow \texttt{DST}(\hat{V}_W^{+}), \quad \Sigma_W \leftarrow \Sigma_W^{+}$ \COMMENT{Transform back via DST.}
\STATE $U_V \leftarrow \texttt{DST}(\hat{U}_V^{+}), \quad V_V \leftarrow \texttt{DST}(\hat{V}_V^{+}), \quad \Sigma_V \leftarrow \Sigma_V^{+}$
\RETURN $[U_W, \Sigma_W, V_W], [U_V, \Sigma_V, V_V]$
\end{algorithmic}
\end{algorithm}

\subsection{The low-rank OFT solver}
\label{sec:lr_oft}
The OFT quadrature \cref{eq:oft_quadrature} requires the weighted sum $\sum_{k=0}^N \omega_k X^k$ of the complex solutions $X^k = V^k + iW^k$ at the nodes $\tau_k$, with complex weights $\omega_k = \omega_k^{\rm re} + i\,\omega_k^{\rm im}$. Since
\begin{equation}
    \label{eq:lr_weighted}
    \omega_k X^k = \left(\omega_k^{\rm re} V^k - \omega_k^{\rm im} W^k\right) + i\left(\omega_k^{\rm im} V^k + \omega_k^{\rm re} W^k\right),
\end{equation}
we store the accumulated sum as $S_V + iS_W$, with $S_V$ and $S_W$ real low-rank matrices in SVD form, and after each time step update
\begin{equation}
    \label{eq:lr_accumulate}
    \begin{aligned}
        S_V &\leftarrow \texttt{truncsum}\left(S_V + \omega_k^{\rm re} V^{k} - \omega_k^{\rm im} W^{k}\right), \\
        S_W &\leftarrow \texttt{truncsum}\left(S_W + \omega_k^{\rm im} V^{k} + \omega_k^{\rm re} W^{k}\right).
    \end{aligned}
\end{equation}
Each argument is a sum of three low-rank matrices in SVD form, where the scalar weights are absorbed into the singular values, so the rank of the sum can be as large as the sum of the three ranks. The rounded sum of \cref{sec:truncsum}, with tolerance $\Delta t_k^3$ (in Example~1 a fixed tolerance of $10^{-12}$ is used instead; see \cref{sec:experiments}), recompresses it, which keeps the rank of $S_V$ and $S_W$ under control throughout the integration.

The full Helmholtz solver is given in \cref{alg:helmholtz_oft}. It evaluates the OFT integral twice, \cref{eq:helmholtz_oft1} followed by \cref{eq:helmholtz_oft2}, with \cref{alg:oft_sqrtinv}. The first pass takes $g$ as initial data and returns $\tilde{v} = f(A)g$ in low-rank form; the second pass takes $\tilde{v}$ as initial data and returns $v = f(A)\tilde{v} = A^{-1}g$. Within each pass, \cref{alg:oft_sqrtinv} calls \cref{alg:timestep} once per quadrature node and updates the accumulated sum with \cref{eq:lr_accumulate}.

\paragraph{Evolution of the rank}
The ranks of $V^k$, $W^k$, $S_V$ and $S_W$ are not fixed in advance. In each step, \texttt{Cross-DEIM} returns the smallest rank for which the discarded singular values of $\hat{V}^{k+1}$ and $\hat{W}^{k+1}$ are below the tolerance $\Delta t_k^3$, and \texttt{truncsum} does the same for $S_V$ and $S_W$. The ranks therefore follow the numerical rank of the computed quantities, and two effects change them over time. First, the factors $\phi_{ij}\cos\theta_{ij}$ and $\phi_{ij}\sin\theta_{ij}$ in \cref{eq:lr_entry_functions} depend on $i$ and $j$ through $\lambda_i + \lambda_j$ and are not products of a function of $i$ and a function of $j$, so a time step can increase the rank of $V^k$ and $W^k$. Second, when $\beta > 0$ we have $\phi_{ij} < 1$, so the singular values of $V^k$ and $W^k$ decay with $\tau$, and the rank decreases once they fall below the tolerance. The rank of the accumulated sums grows while the integrand adds new contributions and levels off once it has decayed. In Example~1 the initial data is a single eigenvector of the discrete Laplacian, so neither effect occurs and the rank remains one (\cref{fig:ex1_ranks}). In Example~2 the ranks first grow and then decrease (\cref{fig:ex2_ranks}).

\paragraph{Computational cost}
Let $r$ denote the largest rank of $V^k$, $W^k$, $S_V$ and $S_W$ during the integration. In each time step, the discrete sine transforms of the factors cost $O(r\, n \log n)$; the two Cross-DEIM calls cost $O(n r^2)$, since they sample $O(nr)$ entries at $O(r)$ each and the cross approximation itself costs $O(nr^2)$; and the two rounded sums, with total rank at most $3r$, cost $O(nr^2)$. The cost per time step is therefore $O(n r \log n + n r^2)$, and the memory required is $O(nr)$. Solving the Helmholtz equation requires two passes of $N$ time steps each, for a total cost of $O(N(n r \log n + n r^2))$. In comparison, carrying out the same time stepping with dense $n \times n$ matrices requires $O(n^2 \log n)$ operations per time step for the sine transforms and $O(n^2)$ memory. The low-rank solver is therefore advantageous when $r \ll n$, that is, when the solution of the Schr\"odinger equation and the accumulated sum remain of low rank throughout the integration. When $r$ grows to be comparable to $n$, this advantage is lost.

\begin{algorithm}[]
\caption{{\tt v = Helmholtz-OFT(g, $\alpha$, $\beta$, $\kappa$, $\Delta t_0$, $\Delta t_f$, $T$)}\\
Low-rank OFT Helmholtz solve: $v = A^{-1}g$ via $v = f(A)^2 g$, where $A = m(\vect{x}) + \Delta/\kappa^2$ and $f(A) = A^{-1/2}$}
\label{alg:helmholtz_oft}
\begin{algorithmic}[1]
\STATE {\bf Input:} Source term $g \in \mathbb{C}^{n \times n}$, parameters $\alpha, \beta, \kappa$, initial step size $\Delta t_0$, final step size $\Delta t_f$, final time $T$
\STATE {\bf Output:} Approximate solution $v \approx A^{-1}g$ in low-rank form $U_V \Sigma_V V_V^T + i\, U_W \Sigma_W V_W^T$
\STATE Compute exponential time grid $\{\tau_k\}$ via \cref{eq:tau_def}
\STATE Compute Fresnel quadrature weights $\{\omega_k\}$ via \cref{eq:OFTWeights}
\STATE $[U_W, \Sigma_W, V_W], [U_V, \Sigma_V, V_V] \leftarrow \texttt{OFT-SqrtInv}(W^{(0)}\text{-factors}, V^{(0)}\text{-factors}, \{\tau_k\}, \{\omega_k\}, \alpha, \beta, \kappa)$ \COMMENT{Pass 1: $v_1 = A^{-1/2}g$.}
\STATE $[U_W, \Sigma_W, V_W], [U_V, \Sigma_V, V_V] \leftarrow \texttt{OFT-SqrtInv}([U_W, \Sigma_W, V_W], [U_V, \Sigma_V, V_V], \{\tau_k\}, \{\omega_k\}, \alpha, \beta, \kappa)$ \COMMENT{Pass 2: $v = A^{-1/2}v_1 = A^{-1}g$.}
\STATE $v \leftarrow U_V \Sigma_V V_V^T + i\, U_W \Sigma_W V_W^T$ \COMMENT{Reconstruct complex solution.}
\RETURN $v$
\end{algorithmic}
\end{algorithm}

\section{Experimental results}
\label{sec:experiments}
Here we will discuss the results of numerical experiments using our method, and compare them to the results obtained when not using low-rank approximation. Errors are reported as relative errors in the Frobenius norm, $\|v - v_{\rm ref}\|_F / \|v_{\rm ref}\|_F$, where $v_{\rm ref}$ is the exact solution sampled on the grid in Example~1 and one of the reference solutions described in \cref{sec:ex2} in Example~2. In both examples the \texttt{Cross-DEIM} tolerance in time step $k$ is $\Delta t_k^3$. The \texttt{truncsum} tolerance is $10^{-12}$ in Example~1 and $\Delta t_k^3$ in Example~2. In all cases the maximum allowed rank is taken to be the full size of the mesh in either dimension.

\subsection{Example 1}
For our first experiment we solve
\begin{equation}
    \label{eq:experiment1}
    u - i \Delta u = f \quad \text{in } \Omega = [0,1]^2, \qquad u = 0 \quad \text{on } \partial\Omega,
\end{equation}
with $f(\vect{x}) = \left(1 + i\pi^2\left(k^2 + l^2\right)\right)\sin\left(k \pi x\right)\sin\left(l \pi y\right)$ and $k = l = 1$. This is the setting of \cref{eq:oft2} with $A = \Delta$, and it corresponds to $\alpha = \beta = 0$ and $\kappa = 1$ in \cref{sec:low-rank}. We choose it because its exact solution,
\begin{equation}
    \label{eq:ex1_exact_sol}
    u(\vect{x}) = \sin\left(k \pi x\right)\sin\left(l \pi y\right),
\end{equation}
is of rank one. By \cref{eq:oft2}, $u = \int_0^\infty e^{-\tau} \hat{u}(\tau)\, d\tau$, where $\hat{u}(\tau) = e^{i\tau\Delta} f$ solves the initial value problem
\begin{equation}
    \label{eq:exp_ivp}
    \hat{u}_\tau = i \Delta \hat{u}, \qquad \hat{u}(0) = f,
\end{equation}
which we evaluate with the low-rank scheme of \cref{sec:low-rank,sec:lowrank_impl}. Unlike Example~2, this example uses uniform steps $\tau_j = j\Delta t$, $j = 0, \dots, N$, with $N\Delta t = T$, and approximates the integral by the trapezoidal rule,
\begin{equation}
    \label{eq:exp_weighted_sum}
    u \approx \sum_{j=0}^{N} w_j\, \hat{u}(\tau_j), \qquad w_j = \Delta t\, e^{-\tau_j} \ \text{for } 0 < j < N, \qquad w_0 = \tfrac{\Delta t}{2}, \quad w_N = \tfrac{\Delta t}{2} e^{-T}.
\end{equation}

\begin{figure}[t]
    \centering
    \begin{minipage}[b]{0.64\textwidth}
        \centering
        \includegraphics[width=0.49\linewidth, trim={58bp 15.5bp 9bp 83bp}, clip]{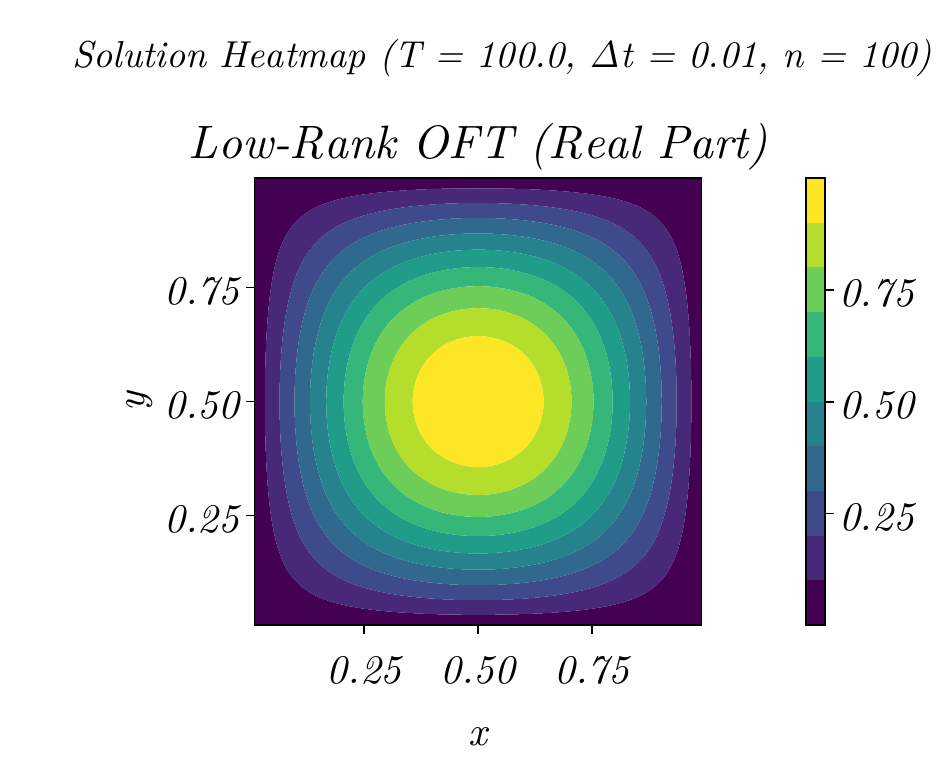}\hfill
        \includegraphics[width=0.49\linewidth, trim={58bp 15.5bp 9bp 83bp}, clip]{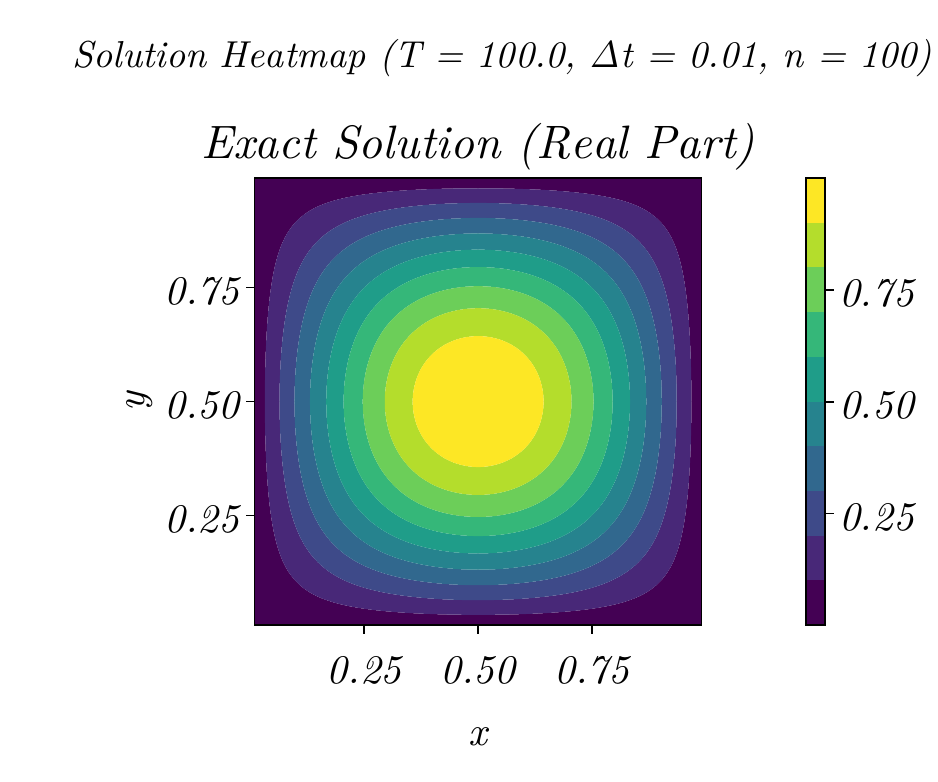}
        \caption{Real part of the computed low-rank solution to \cref{eq:experiment1} (left) and of the exact solution \cref{eq:ex1_exact_sol} (right), with $T=100$, $\Delta t = 0.01$, and $n = 100$.}
        \label{fig:ex1_sol}
    \end{minipage}\hfill
    \begin{minipage}[b]{0.33\textwidth}
        \centering
        \includegraphics[width=\linewidth, trim={17bp 15.5bp 10bp 41.5bp}, clip]{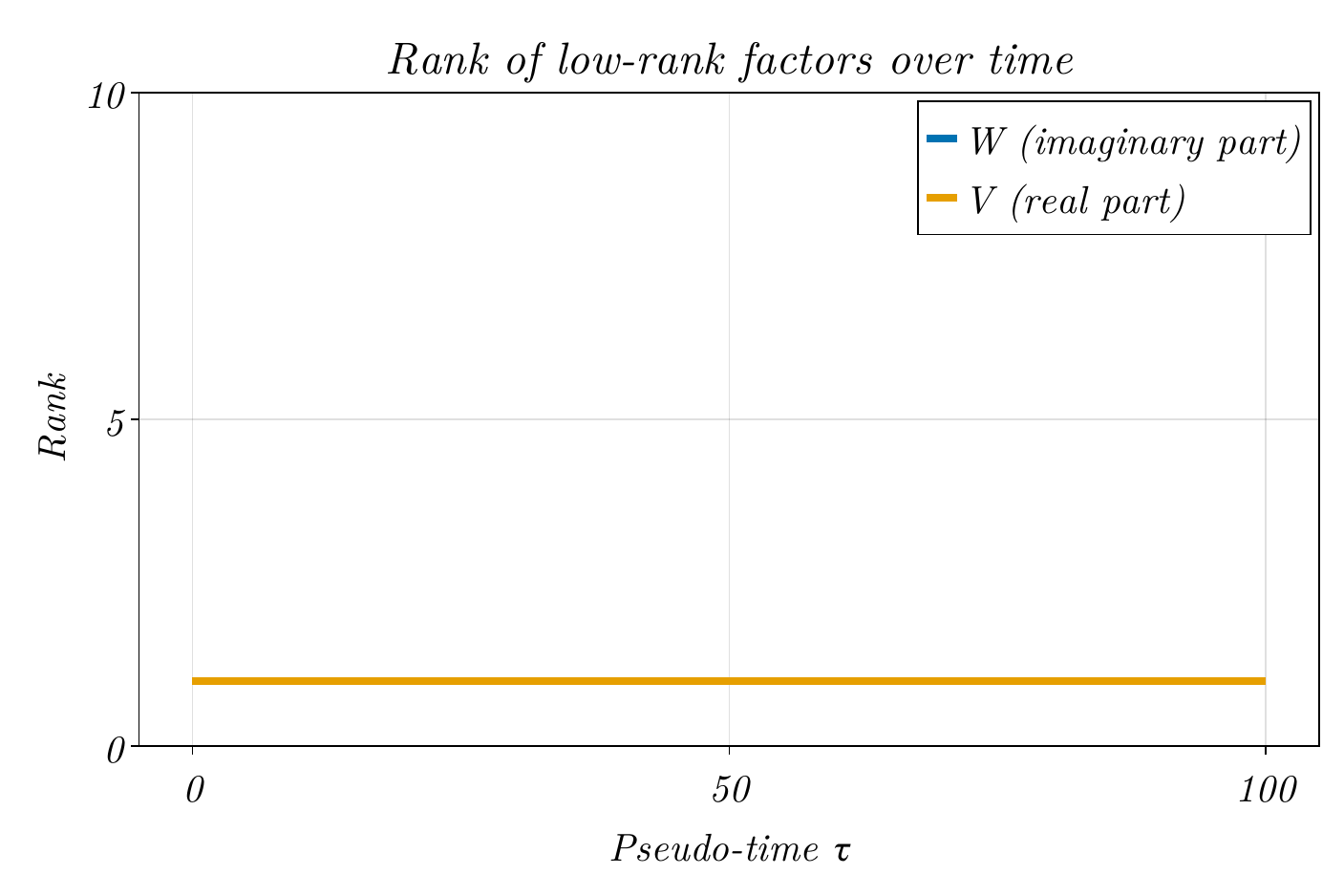}
        \caption{Rank of $V^k$ and $W^k$ in Example~1.}
        \label{fig:ex1_ranks}
    \end{minipage}
\end{figure}

\begin{figure}[t]
    \centering
    \begin{minipage}[t]{0.48\textwidth}
        \centering
        \includegraphics[width=\linewidth, trim={17.5bp 11.5bp 10bp 44.5bp}, clip]{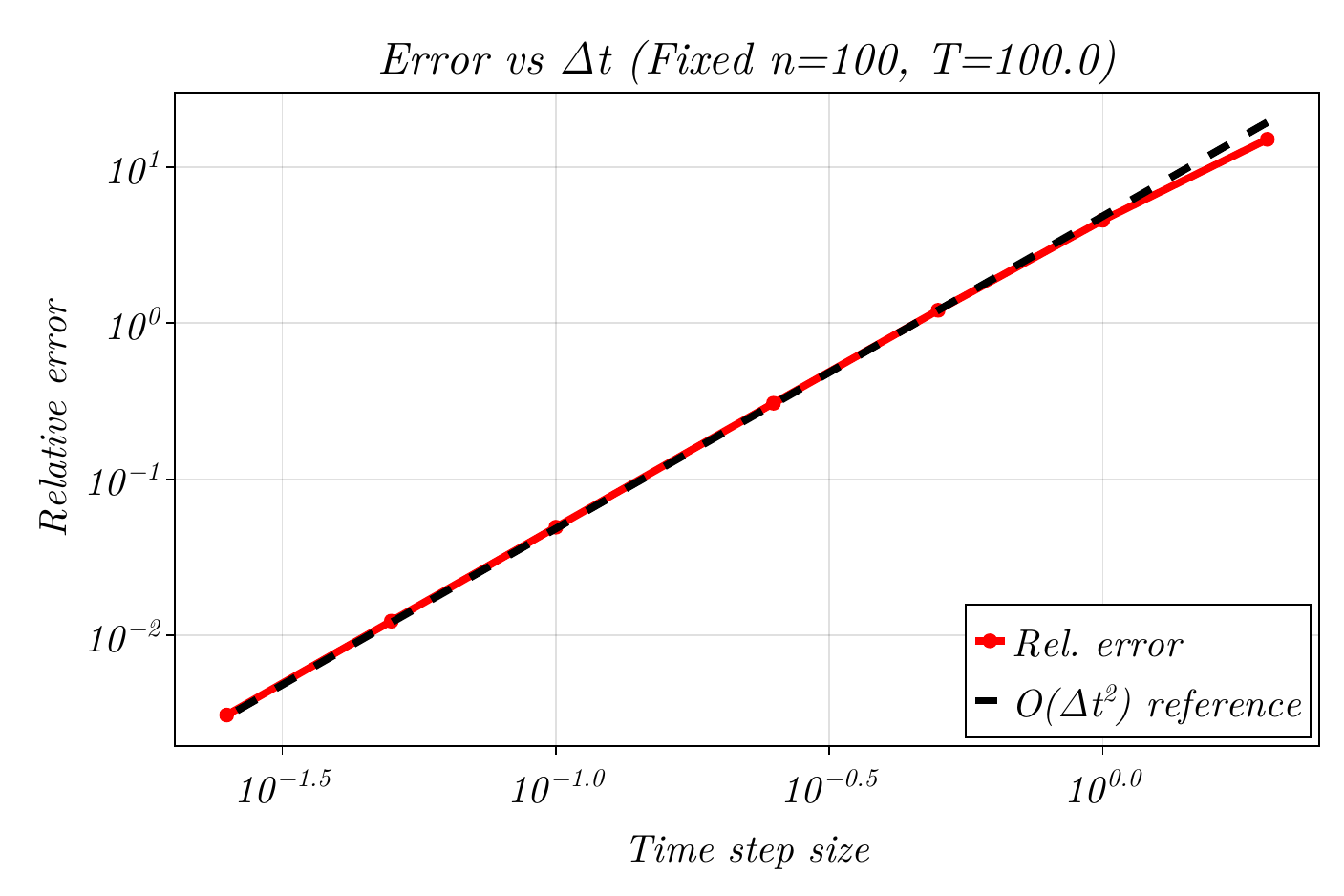}
        \caption{Error as a function of the time step $\Delta t$, with $T=100$ and $n = 100$. The dashed line shows second order.}
        \label{fig:ex1_err_dt}
    \end{minipage}\hfill
    \begin{minipage}[t]{0.48\textwidth}
        \centering
        \includegraphics[width=\linewidth, trim={17bp 11.5bp 10bp 45.5bp}, clip]{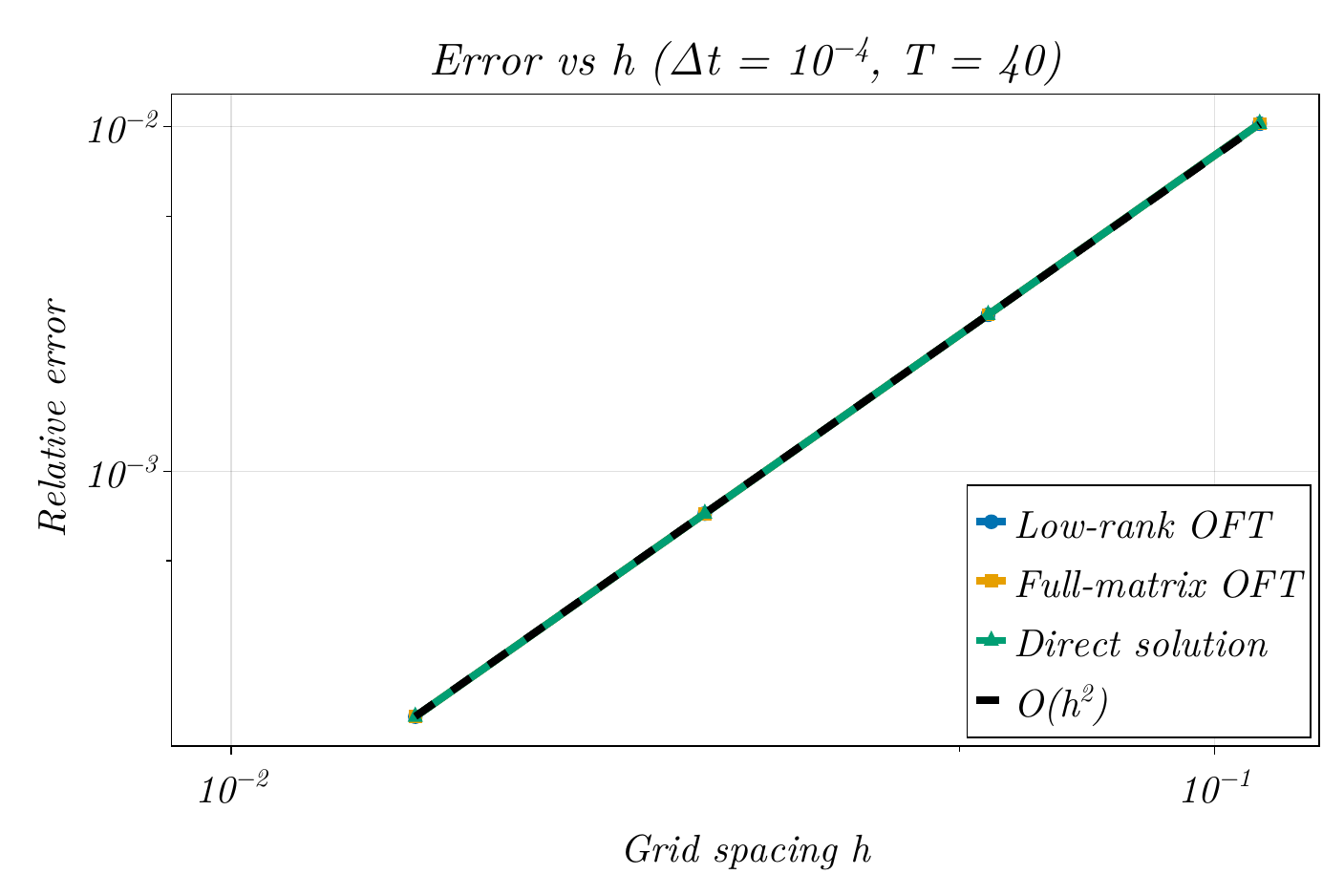}
        \caption{Error as a function of $h$, with $\Delta t = 10^{-4}$ and $T = 40$: low-rank OFT, full-matrix OFT, and direct solution of the discretized linear system. The dashed line shows second order.}
        \label{fig:ex1_err_h}
    \end{minipage}
\end{figure}

\Cref{fig:ex1_sol} compares the computed low-rank solution with the exact solution \cref{eq:ex1_exact_sol}. \Cref{fig:ex1_ranks} shows that the rank remains one over all time steps. This is because $f$, sampled on the grid, is a single eigenvector of the discrete Laplacian, so the solution of the spatially discretized problem \cref{eq:exp_ivp} at each $\tau$ is $f$ multiplied by a complex phase factor $e^{i\tau\mu_{kl}}$, where $\mu_{kl} = -(\lambda_k + \lambda_l)/h^2 < 0$ is the corresponding eigenvalue of the discrete Laplacian, and no new modes are introduced.

\Cref{fig:ex1_err_dt} shows that the error decreases at second order in $\Delta t$, as expected since both the trapezoidal rule and the midpoint rule in pseudo-time are second-order accurate. \Cref{fig:ex1_err_h} shows the error as the grid is refined, with $\Delta t$ small enough that the time discretization error and the truncation of the integral are negligible. The error decreases at second order, and the low-rank OFT, the full-matrix OFT and the direct solution of the linear system obtained from the second-order finite-difference discretization of \cref{eq:experiment1} agree to at least five significant digits. Since $f$ is a single eigenfunction of the discrete Laplacian, the error is dominated by the error in the approximation of the corresponding eigenvalue.

\subsection{Example 2}
\label{sec:ex2}
We now solve the Helmholtz equation \cref{eq:helmholtz} on $\Omega = [-\pi, \pi]^2$, so that $L = 2\pi$ and $h = 2\pi/(n+1)$ in \cref{eq:grid}, with the source $g(\vect{x}) = \exp\left(-6\left(x^2 + y^2\right)\right)$, free-space wavenumber $\kappa = 2$, and the constant coefficient $m = 3 + 0.1\, i$, i.e.\ $\alpha = 2$ and $\beta = 0.1$ in \cref{eq:base_eq_split}. The solution is computed with \cref{alg:helmholtz_oft}, that is, with two passes of the OFT integral, \cref{eq:helmholtz_oft1} followed by \cref{eq:helmholtz_oft2}, each of which solves \cref{eq:base_eq} with initial data $g$ and $\tilde{v}$, respectively. Unlike Example~1, this example uses the exponential time stepping of \cref{sec:num_int_oft}, with $\Delta t_f = 10\,\Delta t_0$ in all runs. Since $|g| < 10^{-25}$ on $\partial\Omega$, the homogeneous Dirichlet condition is consistent with $g$ to machine precision. Unless stated otherwise, all runs use $T = 100$ and $n = 200$.

\paragraph{Reference solutions}
The first reference solution, $v_h$, is the direct solution of the linear system obtained by discretizing \cref{eq:helmholtz} on the same grid with the second-order finite differences \cref{eq:five_point},
\begin{equation}
    \label{eq:ex2_linear_system}
    m\, v_h - \frac{1}{h^2\kappa^2}\left(T_n v_h + v_h T_n\right) = g,
\end{equation}
computed with a sparse direct solver. Since the OFT applies exactly this discrete operator, $v - v_h$ contains only the errors of the time discretization, of the truncation of the integral at $T$, and of the low-rank approximation. The second, $v_{\rm fine}$, is computed with the exact eigenvalues of the Laplacian on a grid with at least 2047 interior points in each direction and restricted to the coarse grid, so $v - v_{\rm fine}$ also contains the spatial discretization error. To isolate the low-rank error, we also compute the full-matrix OFT, which uses the same quadrature, time steps and update rule \cref{eq:theta_update} as the low-rank OFT but applies it to all $n^2$ sine coefficients, without any low-rank truncation.

\begin{figure}[t]
    \centering
    \includegraphics[width=0.4\textwidth, trim={56.5bp 15.5bp 9bp 83bp}, clip]{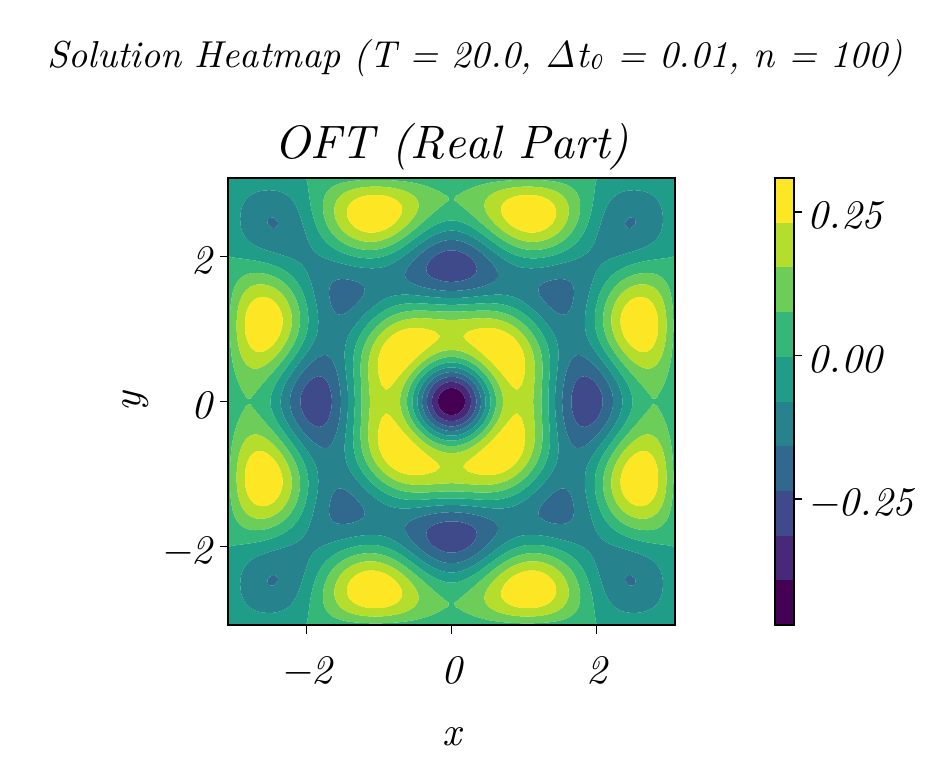}\qquad
    \includegraphics[width=0.4\textwidth, trim={56.5bp 15.5bp 9bp 83bp}, clip]{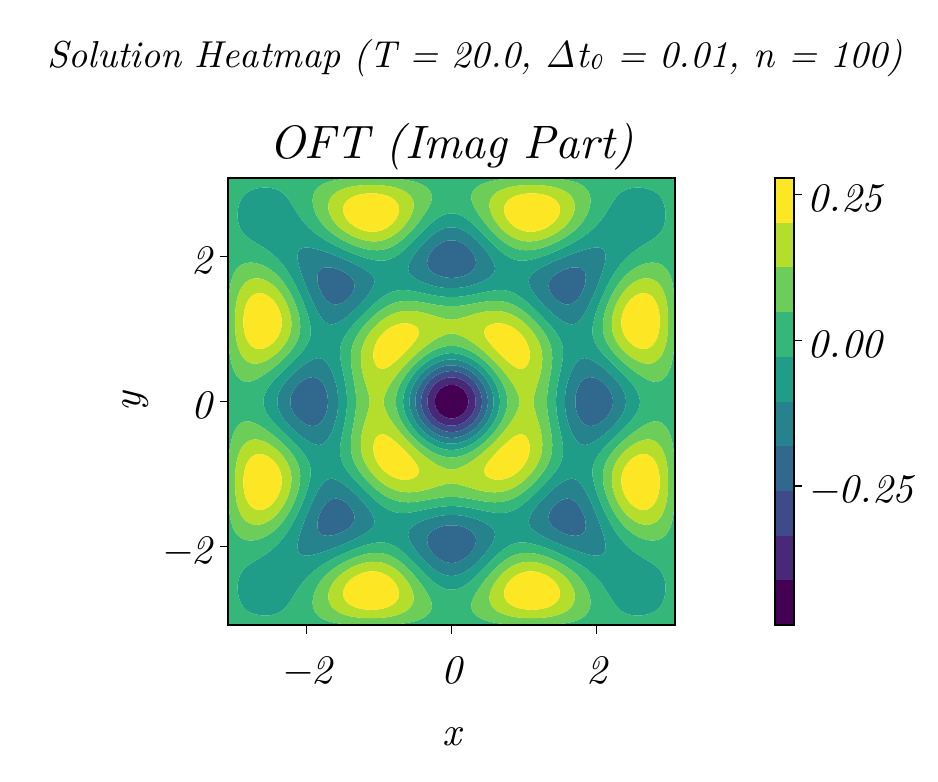}
    \caption{Real (left) and imaginary (right) parts of the computed low-rank solution $v$ in Example~2 with $T = 20$, $\Delta t_0 = 0.01$, $\Delta t_f = 0.1$, and $n = 100$.}
    \label{fig:ex2_sol}
\end{figure}

\begin{figure}[t]
    \centering
    \includegraphics[width=0.9\textwidth, trim={17bp 15.5bp 10bp 67.5bp}, clip]{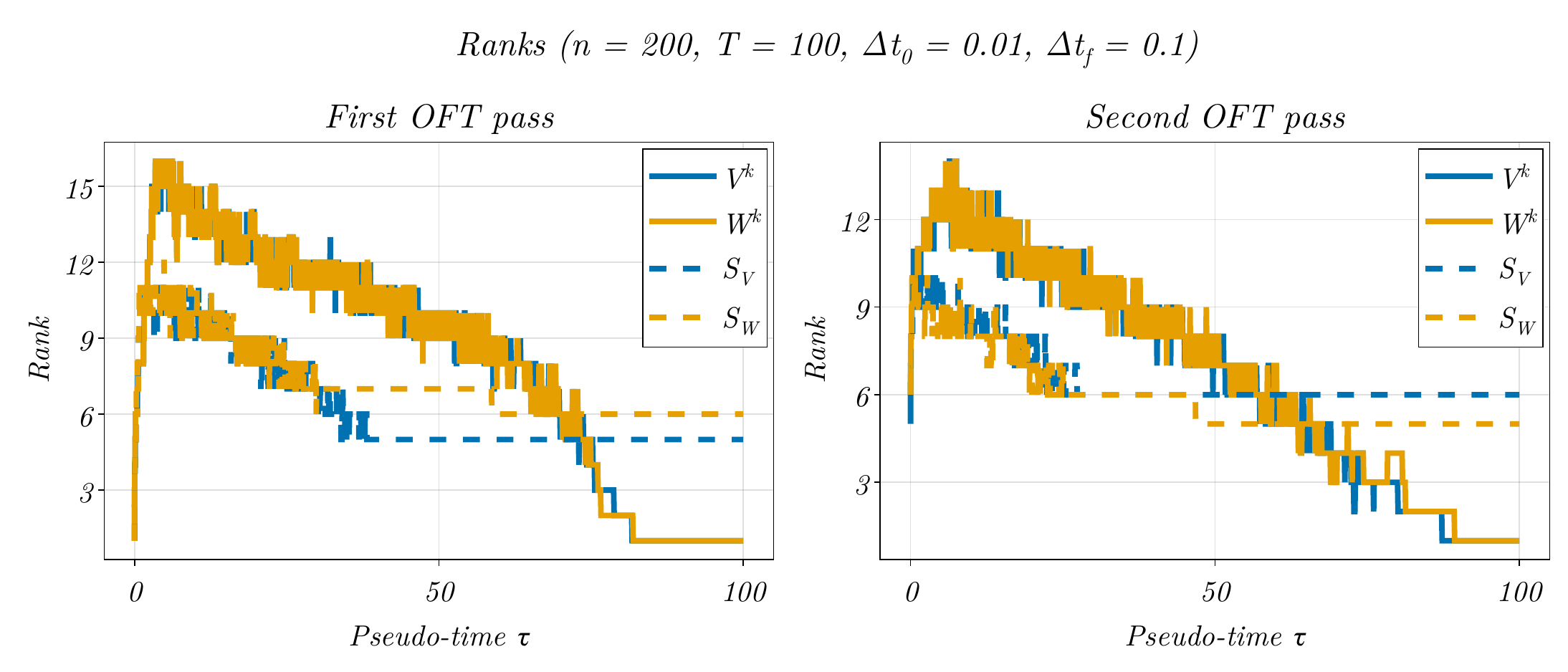}
    \caption{Ranks of $V^k$ and $W^k$ (solid) and of the accumulated sums $S_V$ and $S_W$ (dashed) during the first (left) and second (right) OFT pass, with $n = 200$, $T = 100$, $\Delta t_0 = 0.01$, and $\Delta t_f = 0.1$.}
    \label{fig:ex2_ranks}
\end{figure}

\Cref{fig:ex2_sol} shows the real and imaginary parts of the computed solution.

\paragraph{Ranks}
\Cref{fig:ex2_ranks} shows the ranks during both OFT passes. The Gaussian source is separable, so the first pass starts from rank one, and the second from the rank of the result of the first pass, five or six. With constant coefficients and separable initial data, the exact solution of the Schr\"odinger equation \cref{eq:base_eq} remains separable, so its real and imaginary parts $V$ and $W$ have rank at most two. The computed $V^k$ and $W^k$ reach ranks of up to 16. To investigate this, we repeated the time stepping with full matrices: at $\tau \approx 3.4$, where the rank peaks, the solution computed with \cref{eq:theta_update} has 16 singular values above the tolerance, while the exact solution of the discretized Schr\"odinger equations has two. This suggests that the additional rank comes from the time discretization rather than from the Schr\"odinger equation itself. As $\tau$ increases, the absorption damps the solution and the ranks decrease again. The accumulated sums reach ranks of at most 12 and end with rank five or six, which reflects the rapid decay of the singular values of the Helmholtz solution itself.

\begin{figure}[t]
    \centering
    \begin{minipage}[t]{0.48\textwidth}
        \centering
        \includegraphics[width=\linewidth, trim={17bp 14bp 10bp 48bp}, clip]{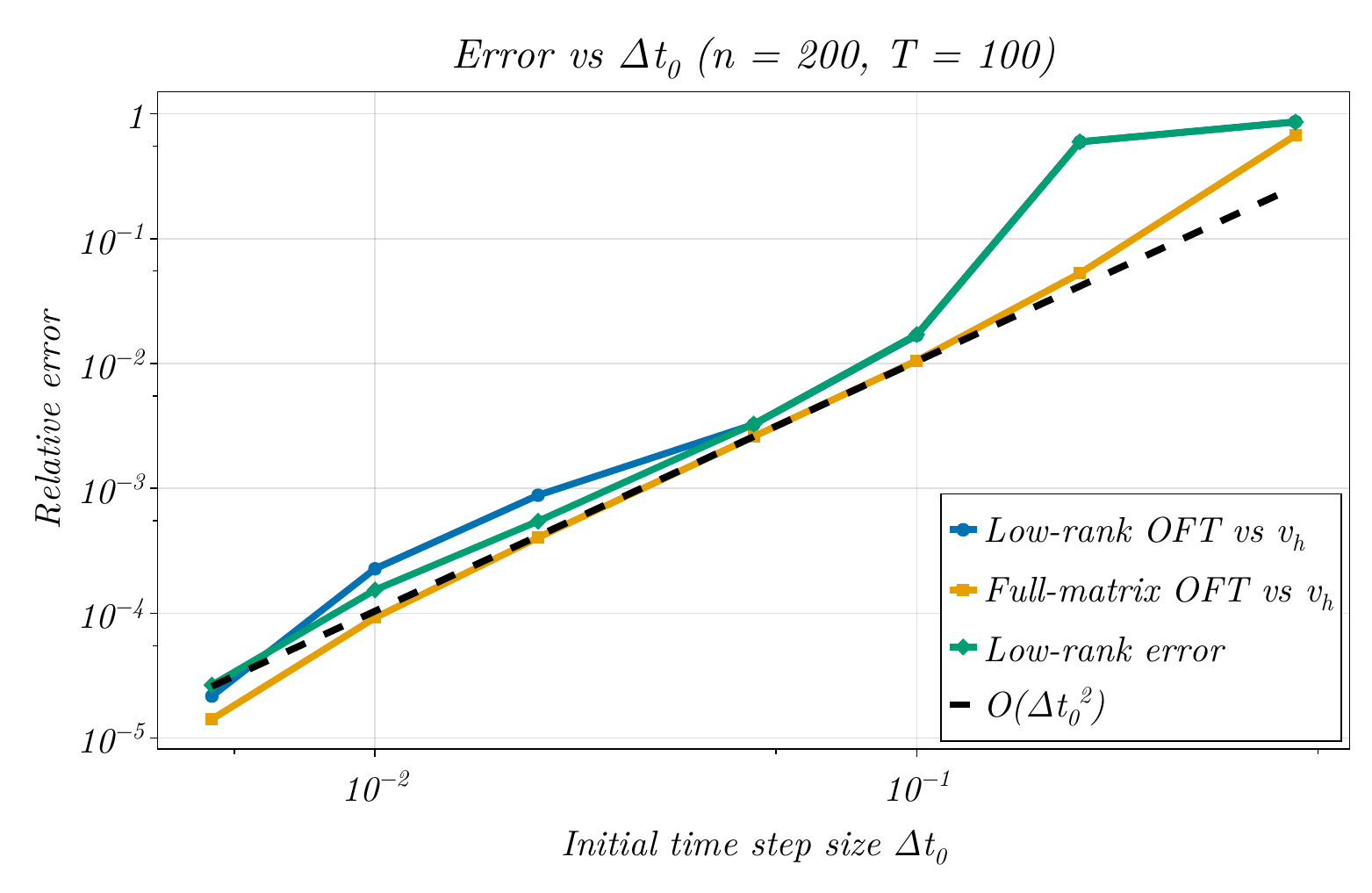}
        \caption{Error as a function of $\Delta t_0$: low-rank and full-matrix OFT measured against $v_h$, and the difference between the two (the low-rank error). The dashed line shows second order.}
        \label{fig:ex2_err_dt}
    \end{minipage}\hfill
    \begin{minipage}[t]{0.48\textwidth}
        \centering
        \includegraphics[width=\linewidth, trim={17.5bp 12bp 10bp 43.5bp}, clip]{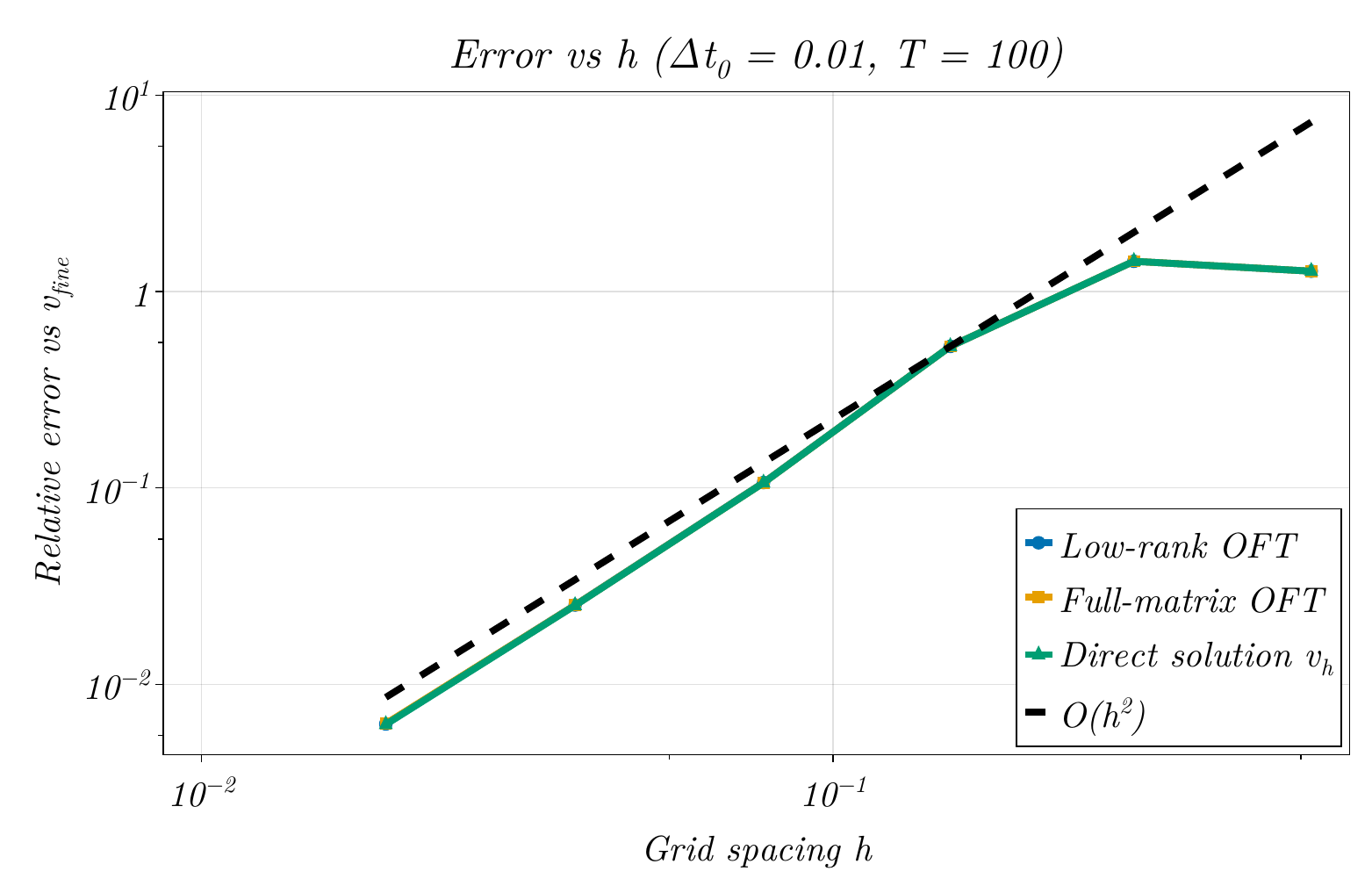}
        \caption{Error as a function of $h$, measured against $v_{\rm fine}$, for $\Delta t_0 = 0.01$: low-rank OFT, full-matrix OFT, and direct solution $v_h$ of \cref{eq:ex2_linear_system}. The dashed line shows second order.}
        \label{fig:ex2_err_h}
    \end{minipage}
\end{figure}

\paragraph{Convergence}
\Cref{fig:ex2_err_dt} shows the error as a function of $\Delta t_0$. Measured against $v_h$, the full-matrix OFT converges at second order, from $1.1 \times 10^{-2}$ at $\Delta t_0 = 0.1$ to $1.4 \times 10^{-5}$ at $\Delta t_0 = 0.005$. The low-rank error decreases at approximately the same rate, and the low-rank OFT reaches an error of $2.2 \times 10^{-5}$ at $\Delta t_0 = 0.005$. For $\Delta t_0 \ge 0.2$ the low-rank approximation fails, with relative errors above 0.5; at $\Delta t_0 = 0.5$ the full-matrix OFT is also inaccurate. The failure is caused by the size of the tolerance $\epsilon_k = \Delta t_k^3$, which is an absolute, not a relative, tolerance. In each time step, \texttt{Cross-DEIM} truncates the approximation of $\hat{V}^{k+1}$ and $\hat{W}^{k+1}$ to the smallest rank $r^\ast$ for which the discarded singular values satisfy $\sum_{l > r^\ast} S_l^2 < \epsilon_k^2$, so it introduces an error of at most $\epsilon_k$ in the Frobenius norm, and \texttt{truncsum} discards all singular values of $S_V$ and $S_W$ that are smaller than $\epsilon_k$. By \cref{eq:tau_def,eq:tau_params}, the step sizes grow linearly in $\tau$,
\begin{equation}
    \label{eq:dt_linear}
    \Delta t_k = \tau_{k+1} - \tau_k = (b - 1)(\tau_k + a) = \Delta t_0 + \left(\Delta t_f - \Delta t_0\right)\frac{\tau_k}{T},
\end{equation}
so in the last steps $\Delta t_k \approx \Delta t_f = 10\,\Delta t_0$ and $\epsilon_k \approx 1000\,\Delta t_0^3$. The accumulated sums approximate $f(A)g$ and $v$, whose leading singular values are of order one to ten; for example, the singular values of $\operatorname{Re} f(A)g$ are $10.3$, $2.3$, $1.9$ and $0.95$, and the remaining ones are below $10^{-2}$. For $\Delta t_0 \le 0.1$ we have $\epsilon_k \le 1$, so at most the smallest of these is affected, and only in the last steps. For $\Delta t_0 = 0.2$, however, $\epsilon_k$ reaches 8 and \texttt{truncsum} removes all but the one or two largest singular values of the accumulated sums, which introduces an error of order one in the solution.

\Cref{fig:ex2_err_h} shows the error measured against $v_{\rm fine}$ as the grid is refined. Here the spatial discretization error dominates. For $n \ge 80$, the error of the low-rank OFT agrees with that of $v_h$ to within 1\% and decreases at second order, from $1.1 \times 10^{-1}$ at $n = 80$ to $6.3 \times 10^{-3}$ at $n = 320$; for $n \le 40$ the grid is too coarse to resolve the solution. In summary, at $n = 200$ and $\Delta t_0 = 0.01$ the total error of $1.6 \times 10^{-2}$ is dominated by the spatial discretization, while the low-rank approximation contributes $1.5 \times 10^{-4}$, and the time discretization together with the truncation of the integral at $T$ contributes $9.2 \times 10^{-5}$ (the error of the full-matrix OFT against $v_h$; see \cref{tab:ex2_tol}). The accuracy of the Helmholtz solution itself is therefore limited by the second-order finite differences: measured against $v_{\rm fine}$, the smallest error in \cref{fig:ex2_err_h} is $6.3 \times 10^{-3}$, about two digits. The smaller errors in \cref{fig:ex2_err_dt} are measured against $v_h$ and show only that the time discretization, the truncation of the integral and the low-rank approximation do not limit the accuracy. Since the spatial error is $O(h^2)$, each additional digit requires $n$ to grow by a factor of about $\sqrt{10}$, together with a correspondingly smaller $\Delta t_0$, unless a higher-order discretization in space is used.

\paragraph{Choice of tolerance}
With the tolerance $\Delta t_k^3$, each step introduces a low-rank error of at most $\Delta t_k^3$, so that, heuristically, the accumulated low-rank error is bounded by $\sum_k \Delta t_k^3 \le \Delta t_f^2\, T$, which is of the same order as the second-order time discretization error. \Cref{tab:ex2_tol} compares this choice with the fixed tolerance $10^{-8}$. With the fixed tolerance the low-rank error is negligible, below $3 \times 10^{-9}$, and the total error equals that of the full-matrix OFT; a tolerance of $10^{-10}$ gives the same errors with ranks larger by one or two. With $\Delta t_k^3$ the total error is at most 2.5 times larger, while the maximum ranks are lower.

\begin{table}[t]
    \centering
    \caption{Error (against $v_h$), and maximum ranks over both passes of the low-rank OFT in Example~2 for two tolerances, with $n = 200$ and $T = 100$.}
    \label{tab:ex2_tol}
    \begin{tabular}{llccc}
        \hline
        $\Delta t_0$ & Tolerance & Error & Rank $V^k, W^k$ & Rank $S_V, S_W$ \\
        \hline
        0.1  & $\Delta t_k^3$ & $1.7 \times 10^{-2}$ & 12 & 9  \\
             & $10^{-8}$      & $1.1 \times 10^{-2}$ & 17 & 16 \\
        0.05 & $\Delta t_k^3$ & $3.2 \times 10^{-3}$ & 14 & 10 \\
             & $10^{-8}$      & $2.6 \times 10^{-3}$ & 18 & 15 \\
        0.02 & $\Delta t_k^3$ & $8.8 \times 10^{-4}$ & 16 & 11 \\
             & $10^{-8}$      & $4.0 \times 10^{-4}$ & 18 & 15 \\
        0.01 & $\Delta t_k^3$ & $2.3 \times 10^{-4}$ & 16 & 12 \\
             & $10^{-8}$      & $9.2 \times 10^{-5}$ & 19 & 15 \\
        \hline
    \end{tabular}
\end{table}

\paragraph{Rank and grid size}
With $\Delta t_0 = 0.01$, the maximum rank of the accumulated sums in the second OFT pass is 10, 11, 11 and 12 for $n = 64$, 128, 256 and 512. The rank $r$ in the cost estimate of \cref{sec:lr_oft} is the largest rank of all four quantities, including $V^k$ and $W^k$, whose ranks reach 16 at $n = 200$ and are reported here only for that grid. If these ranks also remain bounded as the grid is refined, as those of the accumulated sums do, the cost $O(n r \log n + n r^2)$ per time step grows nearly linearly in $n$, compared with $O(n^2 \log n)$ for the full-matrix OFT.

\section{Conclusions}
\label{sec:conclusions}

Through our experiments we see that the OFT framework lends itself well to the use of low-rank methods. In both examples the low-rank solver converges at second order in the time step, as expected since the midpoint rule \cref{eq:step1} is equivalent to the trapezoidal rule, with the right-hand side evaluated at the average of the two time levels. In Example~2, the Helmholtz equation with an absorbing medium, the low-rank OFT reaches the accuracy of the direct solution of the discretized linear system, while the solutions of the Schr\"odinger equation and the accumulated OFT sums remain of rank at most 16 and 12, respectively, on a $200 \times 200$ grid. Choosing the \texttt{Cross-DEIM} and \texttt{truncsum} tolerances as $\Delta t_k^3$ keeps the error of the low-rank approximation of the same order as the time discretization error, with lower ranks than a fixed small tolerance. The maximum rank of the accumulated sums in the second OFT pass remains between 10 and 12 for grids from $n = 64$ to $n = 512$, which suggests that the cost per time step of the low-rank solver grows nearly linearly in the number of grid points per direction, compared with $O(n^2 \log n)$ for the full-matrix OFT.

Two limitations should be noted. First, the overall accuracy is set by the second-order finite-difference discretization in space, about two digits on the grids used here, so higher accuracy requires finer grids or a higher-order discretization. Second, the time step of \cref{sec:discretization} relies on a constant coefficient $m$. Variable coefficients, and a study of how the ranks of the solutions of the Schr\"odinger equation depend on the grid and on the wavenumber, are left for future work.

\section*{Acknowledgments}
This research is supported by Virginia Tech.

\bibliographystyle{siamplain}
\bibliography{ref}

\appendix
\section{Appendix}

\begin{algorithm}[]
  \caption{{\tt [U,S,V] = Cross-DEIM(G,U0,V0,$\epsilon$, $r_{\rm max}$, $\aleph_{\rm max}$, {\tt maxiter} )}\\
  Adaptive Cross-DEIM approximation to  $G \in \mathbb{R}^{m \times n}$}
  \begin{algorithmic}[1]
  \STATE {\bf Input:} Matrix $G \in \mathbb{R}^{m \times n}$, initial rank $r$ guess to the singular vector matrix $U_0 \in \mathbb{R}^{m \times r}, V_0 \in \mathbb{R}^{n \times r}$, tolerance $\epsilon$, maximum output rank $r_{\rm max}$, maximum index set cardinality $\aleph_{\rm max}$, maximum number of iterations {\tt maxiter}. 
  \STATE {\bf Output:} Approximate SVD of $G$, $U\in \mathbb{R}^{m \times r},$ $S \in \mathbb{R}^{r \times r}$, $V\in \mathbb{R}^{n \times r}$.
  \STATE Set $\mathcal{I}_0=\mathcal{J}_0 = \emptyset$.
\FOR{$k = 1, 2, \ldots$, {\tt maxiter}}
\STATE $\mathcal{I}_{k}^\ast = {\tt QDEIM}(U_{k-1})$
\STATE $\mathcal{J}_{k}^\ast = {\tt QDEIM}(V_{k-1})$ \COMMENT{ {\tt QDEIM} can be replaced by {\tt DEIM}}
\STATE $\mathcal{I}_{k} = \mathcal{I}_{k}^\ast \cup \mathcal{I}_{k-1}, \mathcal{J}_{k} \leftarrow \mathcal{J}_{k}^\ast \cup \mathcal{J}_{k-1}$ \COMMENT{Note that the index sets are ordered by {\tt QDEIM}.}
\IF[Make sure that the index set increases by one]{$|\mathcal{I}_{k}| = |\mathcal{I}_{k-1}|$ or $k=1$} 
\STATE $\mathcal{I}_{k} = \mathcal{I}_{k}^\ast \cup \{ i_{\rm rand} \in \complement (\mathcal{I}_{k}^\ast) \}$ \COMMENT{using a random $i_{\rm rand}$ from the complement of $\mathcal{I}_{k}^\ast$.}
\ENDIF
\IF{$|\mathcal{J}_{k}| = |\mathcal{J}_{k-1}|$ or $k=1$} 
\STATE $\mathcal{J}_{k} \leftarrow \mathcal{J}_{k}^\ast \cup \{ j_{\rm rand} \in \complement (\mathcal{J}_{k}^\ast) \}$ 
\ENDIF
\IF{$|\mathcal{I}_k| > \aleph_{\rm max}$}
\STATE $\mathcal{I}_k \leftarrow \mathcal{I}_k(1:\aleph_{\rm max})$ \COMMENT{Keep the $\aleph_{\rm max}$ most important indices.} 
\ENDIF
\IF{$|\mathcal{J}_k| > \aleph_{\rm max}$}
\STATE $\mathcal{J}_k \leftarrow \mathcal{J}_k(1:\aleph_{\rm max})$ 
\ENDIF
\STATE $[U_k,S_k,V_k,r_{\rm C},r_{\rm R}] = {\tt scross}(G, \mathcal{I}_k,\mathcal{J}_k)$
\FOR{$l = 1, 2, \ldots, |\mathcal{I}_k|$}
\IF{$|(r_{\rm R})_{l}| < 10^{-12}$}
\STATE Remove element $l$ from $\mathcal{I}_k$  \COMMENT{Remove redundant rows in $R = G(\mathcal{I}_k,:)$.}
\ENDIF  
\ENDFOR 
\FOR{$l = 1, 2, \ldots, |\mathcal{J}_k|$}
\IF{$|(r_{\rm C})_{l}| < 10^{-12}$}
\STATE Remove element $l$ from $\mathcal{J}_k$ \COMMENT{Remove redundant columns in $C = G(:,\mathcal{J}_k)$.} 
\ENDIF
\ENDFOR 
\STATE $\phi = \| U_{k}S_{k}V_{k}^T- U_{k-1}S_{k-1}V_{k-1}^T \|, \ \ S_{\rm min} =  \min({\sf diag}(S_k))$
\STATE  $\eta_1 = \| (I(:,\mathcal{I}_k))^T U_k \|_2^{-1}, \ \ \eta_2 = \| V^T_{k} I(\mathcal{J}_k,:) \|_2^{-1}$ 
\IF{$\max(\phi,\min(\eta_1(1+\eta_2),\eta_2(1+\eta_1))S_{\rm min})  < \epsilon$} 
\STATE Break out of for loop  \COMMENT{Above $S_{\rm min}$ is the smallest s.v. in the $k$th approx.}
\ENDIF
\ENDFOR
\STATE Find $r^\ast$ so that $\sum_{l = r^\ast+1}^{\min(m,n)} S_l^2 < \epsilon^2$ 
\STATE Set $r = \max(\min(r^\ast,r_{\rm max}),1)$
\RETURN $U_k(:,1:r), S_k(1:r,1:r), V_k(:,1:r)$
\end{algorithmic}
\label{algo:crossDEIM}
\end{algorithm}

This appendix contains \cref{algo:crossDEIM,alg:oft_sqrtinv}. The Helmholtz solver \cref{alg:helmholtz_oft} calls \cref{alg:oft_sqrtinv} twice to evaluate $f(A) = A^{-1/2}$, which in turn calls \cref{alg:timestep} at each quadrature node to advance the Schr\"odinger equation. The \texttt{truncsum} routine used in \cref{alg:oft_sqrtinv} is described in \cref{sec:truncsum}, and the \texttt{Cross-DEIM} routine used in \cref{alg:timestep} is described in \cref{sec:cross-deim}.

\begin{algorithm}[]
\caption{{\tt [W-factors, V-factors] = OFT-SqrtInv(W$_{\rm in}$-factors, V$_{\rm in}$-factors, $\{\tau_k\}$, $\{\omega_k\}$, $\alpha$, $\beta$, $\kappa$)}\\
Low-rank evaluation of $A^{-1/2}g$ via the OFT integral}
\label{alg:oft_sqrtinv}
\begin{algorithmic}[1]
\STATE {\bf Input:} Initial low-rank factors $U_W^{(0)}\Sigma_W^{(0)}(V_W^{(0)})^T$, $U_V^{(0)}\Sigma_V^{(0)}(V_V^{(0)})^T$; time grid $\{\tau_k\}_{k=0}^N$; quadrature weights $\{\omega_k\}_{k=0}^N$; parameters $\alpha, \beta, \kappa$
\STATE {\bf Output:} Accumulated low-rank factors for $W$ and $V$ components of $A^{-1/2}g$
\STATE $\omega^{\rm re} + i\,\omega^{\rm im} \leftarrow \omega_0, \quad \epsilon_0 \leftarrow (\tau_1 - \tau_0)^3$ \COMMENT{Initialize accumulator with first quadrature point.}
\STATE $[U_V^{\rm acc}, \Sigma_V^{\rm acc}, V_V^{\rm acc}] \leftarrow \texttt{truncsum}\!\left(U_V^{(0)}(\omega^{\rm re} \Sigma_V^{(0)})(V_V^{(0)})^T + U_W^{(0)}(-\omega^{\rm im} \Sigma_W^{(0)})(V_W^{(0)})^T, \epsilon_0, n\right)$
\STATE $[U_W^{\rm acc}, \Sigma_W^{\rm acc}, V_W^{\rm acc}] \leftarrow \texttt{truncsum}\!\left(U_V^{(0)}(\omega^{\rm im} \Sigma_V^{(0)})(V_V^{(0)})^T + U_W^{(0)}(\omega^{\rm re} \Sigma_W^{(0)})(V_W^{(0)})^T, \epsilon_0, n\right)$
\FOR{$k = 0, 1, \ldots, N-1$}
    \STATE $\Delta t_k \leftarrow \tau_{k+1} - \tau_k, \quad \epsilon_k \leftarrow \Delta t_k^3$
    \STATE $[U_W^{(k+1)}, \Sigma_W^{(k+1)}, V_W^{(k+1)}], [U_V^{(k+1)}, \Sigma_V^{(k+1)}, V_V^{(k+1)}] \leftarrow$ \\
    $\qquad \texttt{LR-Timestep}\!\left([U_W^{(k)}, \Sigma_W^{(k)}, V_W^{(k)}], [U_V^{(k)}, \Sigma_V^{(k)}, V_V^{(k)}], \{\lambda_j\}, \Delta t_k, \alpha, \beta, \kappa\right)$ \COMMENT{Advance Schr\"odinger equation one step (\cref{alg:timestep}).}
    \STATE $\omega^{\rm re} + i\,\omega^{\rm im} \leftarrow \omega_{k+1}$ \COMMENT{Accumulate with quadrature weight $\omega_{k+1}$.}
    \STATE $[U_V^{\rm acc}, \Sigma_V^{\rm acc}, V_V^{\rm acc}] \leftarrow \texttt{truncsum}\Big(U_V^{\rm acc}\Sigma_V^{\rm acc}(V_V^{\rm acc})^T + U_V^{(k+1)}(\omega^{\rm re} \Sigma_V^{(k+1)})(V_V^{(k+1)})^T$ \\
    $\qquad +\, U_W^{(k+1)}(-\omega^{\rm im} \Sigma_W^{(k+1)})(V_W^{(k+1)})^T, \epsilon_k, n\Big)$
    \STATE $[U_W^{\rm acc}, \Sigma_W^{\rm acc}, V_W^{\rm acc}] \leftarrow \texttt{truncsum}\Big(U_W^{\rm acc}\Sigma_W^{\rm acc}(V_W^{\rm acc})^T + U_V^{(k+1)}(\omega^{\rm im} \Sigma_V^{(k+1)})(V_V^{(k+1)})^T$ \\
    $\qquad +\, U_W^{(k+1)}(\omega^{\rm re} \Sigma_W^{(k+1)})(V_W^{(k+1)})^T, \epsilon_k, n\Big)$
\ENDFOR
\RETURN $[U_W^{\rm acc}, \Sigma_W^{\rm acc}, V_W^{\rm acc}], [U_V^{\rm acc}, \Sigma_V^{\rm acc}, V_V^{\rm acc}]$
\end{algorithmic}
\end{algorithm}

\end{document}